\documentclass[11pt,a4paper,reqno]{amsart}
\usepackage{graphicx} 
\usepackage{times} 
\usepackage{amsmath} 
\usepackage{amssymb}  
\usepackage{amsthm, dsfont}
\usepackage{latexsym}
\usepackage{amsfonts,bbm}
\usepackage{xcolor}
\usepackage{mathtools}
\usepackage{enumerate}
\usepackage{cite}
\usepackage{tikz}
\usepackage{graphicx,caption}
\usepackage{todonotes}
\usepackage{float}
\usepackage{nicefrac}
\usepackage[foot]{amsaddr}
\usepackage{natbib}

\usepackage{color}
\usepackage{graphicx}

\usepackage{tikz}
\usetikzlibrary{shapes.geometric}

\usepackage{tabularx}
\usepackage{comment}
\usepackage[hidelinks, colorlinks=false]{hyperref}
\usepackage{cleveref}
\usepackage[normalem]{ulem}
\usepackage{algorithm}

\allowdisplaybreaks

\DeclareMathOperator*{\argmin}{arg\,min}

\newcommand{\ceil}[1]{\lceil #1 \rceil}

\usepackage{amsmath,amssymb,amsfonts}
\usepackage{latexsym}
\usepackage{mathrsfs}
\usepackage{bbm}
\usepackage{xcolor}
\usepackage{xcolor}
\usepackage{tikz}
\usepackage{url}
\usepackage{stmaryrd}
\usepackage{dsfont}
\usepackage{placeins}
\usepackage{lipsum}
\usepackage{amsfonts}
\usepackage{graphicx}
\usepackage{epstopdf}
\usepackage{multirow}

\usepackage{amsopn}
\ifpdf
  \DeclareGraphicsExtensions{.eps,.pdf,.png,.jpg}
\else
  \DeclareGraphicsExtensions{.eps}
\fi

\newcommand{\R}{\ensuremath{\mathbb R}}  
\newcommand{\N}{\ensuremath{\mathbb N}}

\newtheorem{remark}{Remark}
\newtheorem{definition}{Definition}
\newtheorem{theorem}{Theorem}
\newtheorem{proposition}{Proposition}
\newtheorem{lemma}{Lemma}
\newtheorem{corollary}{Corollary}
\newtheorem{assumption}{Assumption}

\renewcommand{\k}{\mathsf{k}}
\newcommand{\F}{f^\varepsilon}
\newcommand{\g}{g^\varepsilon}
\newcommand{\G}{G^\varepsilon}

\newcommand{\V}{V^\varepsilon}
\newcommand{\J}{J^\varepsilon}
\newcommand{\B}{B^\varepsilon}
\newcommand{\LF}{L_{f^\varepsilon}}
\newcommand{\rx}{r_{\mathcal{X}}}

\allowdisplaybreaks

\title[Asymptotic Stability of Data-driven Model Predictive Control in the Koopman framework]{Data-driven Model Predictive Control: Asymptotic Stability despite Approximation Errors exemplified in the Koopman framework}
\author{Irene Schimperna$^{1}$}\address{$^{1}$Civil Engineering and Architecture Department, University of Pavia, Italy \\ Mail: \textsc{irene.schimperna01@universitadipavia.it, lalo.magni@unipv.it}}
\author{Karl Worthmann$^{2}$}\address{$^{2}$Optimization-based Control Group, Institute of Mathematics, Technische Universität Ilmenau, Germany.\\ Mail: \textsc{\{karl.worthmann, lea.bold\}@tu-ilmenau.de}}
\author{Manuel Schaller$^{3}$}\address{$^{3}$Faculty of Mathematics, Chemnitz University of Technology\\ Mail: \textsc{manuel.schaller@math.tu-chemnitz.de}}
\author{Lea Bold$^{2}$}
\author{Lalo Magni$^1$}

\thanks{Irene Schimperna would like to thank Karl Worthmann and TU Ilmenau for the hospitality in the period in which this paper was developed. \\
\indent Lea Bold and Karl Worthmann gratefully acknowledge funding by the German Research Foundation (DFG; project numbers $545246093$ and 535860958). }

\begin{document}

\begin{abstract}               
        In this paper, we analyze stability of nonlinear model predictive control (MPC) using data-driven surrogate models in the optimization step. First, we establish asymptotic stability of the origin, a controlled steady state, w.r.t.\ the MPC closed loop without stabilizing terminal conditions for sufficiently long prediction horizons. To this end, we prove that cost controllability of the original system is preserved if sufficiently accurate proportional bounds on the approximation error hold. Here, proportional refers to state and control. The proportionality of the error bounds is a key element to derive asymptotic stability in presence of modeling errors and not only practical asymptotic stability. Second, we exemplarily verify the imposed assumptions for data-driven surrogates generated with kernel extended dynamic mode decomposition based on Koopman operator theory. Hereby, we do not impose invariance assumptions on finite dictionaries, but rather derive all conditions under non-restrictive conditions. Finally, we demonstrate our findings with numerical simulations.
\end{abstract}

\maketitle

\section{Introduction}

\noindent Model predictive control (MPC) is an optimization-based control algorithm, that optimizes the future trajectory of the system on a finite horizon and defines a state feedback law by applying to the system the first element of the computed optimal input sequence~\cite{RawlMayn17}. 
When, e.g., a first-principles-based description of the system under control is not available, data-driven methods can be used to obtain a surrogate model for the prediction and optimization step. 
For linear systems, a popular data-driven approach~\cite{coulson2019data,BerbKohl20} is based on Willems' fundamental lemma~\cite{WillRapi05} avoiding the model identification step and directly using input-output data~\cite{FaulOu23}. 
For nonlinear systems, data-driven surrogate models can be obtained, e.g., by means of nonlinear ARX models \cite{denicolao1997stabilizing}, neural networks \cite{schimperna2024robust-lstm,ren2022tutorial}, Bayesian identification~\cite{piga2019performance} and Gaussian process regression~\cite{hewing2019cautious}.
A numerically and analytically appealing approach is extended dynamic mode decomposition (EDMD~\cite{williams2015data}) based on the Koopman operator. 
The Koopman operator is a linear but infinite dimensional operator that encodes the behavior of the associated nonlinear dynamical system~\cite{mezic2005spectral,rowley2009spectral}. 
EDMD provides a data-driven approximation of the compression of the Koopman operator to a finite-dimensional subspace spanned by the so-called observable functions, see~\cite{Mezi21} and the recently published overview~\cite[Section~2]{StraWort26}. 
Various approaches extend this idea to control systems, see, e.g., the recent preprint~\cite{HaseMezi25} as well as~\cite[Section~3]{StraWort26}. 
In linear EDMD with control (EDMDc; \cite{proctor2016dynamic,korda2018linear}), a surrogate control system, that is linear in the (lifted) state and control, is used. 
However, as shown in~\cite{iacob:toth:schoukens:2022}, this method is often insufficient to capture direct state-control couplings. 
Bilinear models as in~\cite{peitz:otto:rowley:2020,StraScha26} offer a reliable alternative that is underpinned by finite-data error bounds~\cite{nuske2023finite,schaller2023towards}.
Besides ensuring a sufficient amount of data, a key challenge in EDMD is the choice of the observable functions.
Here, data-driven approaches, such as deep neural networks~\cite{yeung2019learning,JohnBala25} or kernel EDMD~\cite{klus2020kernel} provide a remedy.
In kernel EDMD, the dictionary consists of the canonical features of an a-priori chosen kernel function centered at the data points. For this method, it is possible to build upon the solid approximation-theoretic foundation of reproducing kernel Hilbert spaces (RKHSs) to obtain finite-data error bounds~\cite{kohne2024infty} leveraging Koopman invariance of suitably chosen RKHSs. 
This method and its error bounds have been extended to control systems in~\cite{BoldPhil24}. 
An alternative approach, that shares the advantage of flexible state-control sampling, was proposed in~\cite{bevanda2024nonparametric} using an additional kernel function to express the dependency on the control input. So far, however, no error bounds are available for this method.
Koopman operator-based methods are not the only class of surrogate models for which error bounds can be provided. Other appealing identification technique, that could be used to derive data-driven surrogate models with error bounds include kernel and Gaussian-process regression, see~\cite{care2023kernel, bisiacco2024learning} and \cite{lederer2019uniform, scampicchio2025gaussian}, respectively.

The use of data-driven models gives rise to approximation errors, that are related to the choice of the model class used for approximation and the finite amount of data available for the identification. 
While asymptotic stabilization of a steady-state equilibrium can be guaranteed in the nominal case, i.e., when the model used in the optimization step coincides with the system under control, in presence of modeling errors typically only practical asymptotic stability is achieved.
Practical asymptotic stability means that the tracking error decreases until a certain threshold and, then, stagnates. 
Typically, the stagnation can be related proportionally to the approximation error. 
Practical asymptotic stability of MPC with surrogate models based on EDMD has been shown both without and with the use of terminal ingredients in~\cite{bold2024data} and~\cite{WortStra24}, respectively. 
A method to obtain offset-free tracking in presence of modeling errors and asymptotically constant disturbances is so-called offset-free MPC \cite{pannocchia2003disturbance}. 
This technique has been successfully within the Koopman framework. 
For instance, \cite{chen2022offset} considers linear offset-free MPC with surrogate models based on EDMDc, while \cite{SchiBold25} considers bilinear EDMD models and addresses the case in which a full description of the reference steady state is not available.
Finally, asymptotic stability of the MPC closed-loop can be proven under suitable conditions on the modeling error. 
In particular, \cite{santos2008tool} proves asymptotic stability of an MPC closed-loop with soft constraints and terminal conditions for sufficiently small model-error bounds under the assumption that the terminal constraint is met by the MPC using the model. 
Similarly, the recent work \cite{KuntRawl25} derives asymptotic stability of MPC with terminal conditions and without state constraints in presence of a sufficiently small parametric error in the model such that the equilibrium in the origin is preserved, see also~\cite{schimperna2025stability} for the treatment of state constraints and a rigorous verification of all assumptions using Koopman operator theory. However, asymptotic stability of MPC without terminal conditions in the presence of approximation error has not been studied before.

The first contribution of this work is a general framework for stability analysis of MPC with data-driven surrogate models subject to approximation errors.
To be more precise, we show that, if a proportional error bound on the surrogate holds, the MPC algorithm designed using the data-driven surrogate model asymptotically stabilizes the system. With \textit{proportional} we mean that the error vanishes at the desired set point and its magnitude is bounded by a quantity depending on the size of state and input. 
We show that proportional error bounds imply the preservation of cost controllability~\cite{GrimMess05,GrunPann10,CoroGrun20} in data-driven approximations provided a sufficiently high approximation accuracy.
Then, invoking this controllability property, we rigorously show asymptotic stability of the origin w.r.t.\ the MPC closed loop using the surrogate model in the optimization step if the prediction horizon is sufficiently long, where the origin serves as a prototypical controlled steady state.
The second main contribution of this paper is to verify the assumptions on the data-driven surrogate leveraging recently proposed finite-data error bounds on kernel EDMD under non-restrictive assumptions.
To this end, we consider the learning framework proposed in~\cite{BoldPhil24}, suitably modified to ensure that the model is consistent with the system dynamics in the origin. 
We derive proportional and uniform error bounds --~without requiring invariance assumptions on finite dictionaries using a flexible sampling technique to alleviate data requirements. 
Finally, our findings are verified with numerical simulations.\footnote{We point out that this part is an extension of our conference paper~\cite{BoldScha25} now properly embedded into a more general learning framework such that we were able to infer asymptotic stability (instead of only practical asymptotic stability).} We point out that the stability framework for data-driven MPC is independent of the Koopman operator such that it can also be applied with different learning-based techniques.

The paper is organized as follows. In Section~\ref{sec:problem_formulation}, we formulate the data-driven model, the bounds on the approximation error and the MPC algorithm. Asymptotic stability of MPC is proved in Section~\ref{sec:MPC_stability}. 
In Section~\ref{sec:kEDMD}, we formulate the kEDMD surrogate model and show that it respects the required error bounds. 
The proposed algorithm is illustrated in Section~\ref{sec:simulations} with numerical examples. 
Finally, in Section~\ref{sec:conclusions}, we draw the conclusions of the work.

\noindent \textbf{Notation}. $\|\cdot\|$ denotes the Euclidean norm on~$\R^n$ and its induced matrix norm on $\R^{n\times n}$, while $\|\cdot\|_F$ is used for the Frobenius norm. Given two numbers~$a,b$, the abbreviation $[a:b] := \mathbb{Z} \cap [a,b]$ is used. 
For a set $\Omega$, $\overline{\Omega}$ and $\operatorname{int}(\Omega)$ denote its closure and its interior, resp.
For radius~$r > 0$, $\mathcal{B}_r(x_0) := \{ x \in \R^n : \|x - x_0\| \leq r \}$ is the closed ball
centered in~$x_0 \in \R^n$. 
$\mathcal{C}_b(\Omega) \coloneqq \mathcal{C}_b(\Omega; \R^n)$ denotes the space of continuous and bounded functions on~$\Omega$. 
$J_{f}$ denotes the Jacobian of the function $f: \mathbb{R}^n \rightarrow \mathbb{R}$, and $J_{f,x}$ is used to denote the Jacobian of $f: \mathbb{R}^n \times \mathbb{R}^m \rightarrow \mathbb{R}$, $x \mapsto f(x,u)$ w.r.t.\ the argument~$x$.
With $\lambda_{\max}(Q)$ and $\lambda_{\min}(Q)$ we denote the maximum and minimum eigenvalue of the positive definite matrix $Q$, respectively.
The symbols $\oplus$ and $\ominus$ represent the Pontryagin set sum and difference, respectively.

\section{Data-driven MPC formulation}
\label{sec:problem_formulation}

\noindent Consider a discrete-time nonlinear control system given by
\begin{align}\label{eq:system}
    x^+ = f(x, u) 
\end{align}
with state $x \in \R^n$ and input $u \in \mathbb{U}$, where $\mathbb{U} \subseteq \R^m$ is a convex and compact set containing the origin, i.e., $0 \in \R^m$, in its interior.
The system dynamics~$f$ is assumed to be continuous and locally Lipschitz continuous in its first argument, i.e., for each compact set~$K \subset \mathbb{R}^n$ there is a Lipschitz constant $L_f = L_f(K)$ such that, for all $x,y \in K$,
\begin{equation}\label{eq:dynamics:Lipschitz}
    \| f(x,u) - f(y,u) \| \leq L_f \| x - y \| \qquad\forall\, u \in \mathbb{U}.
\end{equation}
Let the origin be a controlled equilibrium of system~\eqref{eq:system} for control input $u = 0$, i.e. $f(0, 0) = 0$.
The objective is to steer the system to the origin using an MPC controller while taking input constraints $u \in \mathbb{U}$ into account.\footnote{The upcoming analysis can be directly applied to the regulation of arbitrary controlled equilibria by considering the shifted dynamics.} However, since the system dynamics~\eqref{eq:system} is assumed to be unknown, the MPC design is based on data-driven surrogate models
\begin{align}\label{eq:model}
    x^+ = f^\varepsilon(x, u)
\end{align}
with $\varepsilon \in (0,\bar{\varepsilon}]$ for some $\bar{\varepsilon} > 0$.
The superscript~$\varepsilon$ refers to the approximation accuracy and, throughout the paper, will be used to refer to all the parameters related to the surrogate dynamics~\eqref{eq:model}. 
Typically, the approximation error depends on the available data (quantity, distribution etc.). 

We consider models accompanied by proportional error bounds, where \textit{proportional} means that the error is zero at the origin and grows, at most, proportionally in the size of state and input. 
This property is formalized in the following assumption.
\begin{assumption}
[Error bounds on a set~$\Omega$]\label{ass:error_bound}
    Consider a set $\Omega \subseteq \mathbb{R}^n$ containing the origin $0_{\mathbb{R}^n}$ in its interior.
    Let $f^\varepsilon$ defined by~\eqref{eq:model} be a surrogate model for the control system~\eqref{eq:system}. 
    For every $\varepsilon \in (0,\bar{\varepsilon}]$, let the surrogate model satisfy the\\
    (\textbf{P}) \textit{proportional} error bound \begin{equation}\label{eq:ass:error_bound:proportional}\tag{P-bound}
                \|f(x, u) - f^\varepsilon(x, u)\| \leq c_x^\varepsilon\|x\| + c_u^\varepsilon\|u\| 
            \end{equation}    
    (\textbf{U}) \textit{uniform} error bound 
    \begin{equation}\label{eq:ass:error_bound:uniform}\tag{U-bound}
        \|f(x, u) - f^\varepsilon(x, u)\| \leq \eta^\varepsilon
    \end{equation}           
    for all $x \in \Omega$, $u \in \mathbb{U}$ such that the parameters $c_x^\varepsilon$, $c_u^\varepsilon$, and $\eta^\varepsilon$ satisfy $\lim_{\varepsilon \searrow 0} \max \{c_x^\varepsilon, c_u^\varepsilon, \eta^\varepsilon \} = 0$.
\end{assumption}
\noindent The proportional bound \eqref{eq:ass:error_bound:proportional} of Assumption~\ref{ass:error_bound} requires that the surrogate model is exact at the origin, which is a reasonable assumption since the control objective is to stabilize the equilibrium located at the origin. 
While a proportional error bound is highly beneficial close to the origin, it may become rather conservative if the distance to the desired set point, i.e., the equilibrium located at the origin, grows. 
Then, using the uniform bounds \eqref{eq:ass:error_bound:uniform} of Assumption~\ref{ass:error_bound} on the approximation error may result in tighter bounds. 

In addition to Assumption~\ref{ass:error_bound}, we require uniform Lipschitz continuity of the surrogate model in the approximation accuracy~$\varepsilon$ to derive asymptotic stability of the origin w.r.t.\ the MPC closed loop.
\begin{assumption}[Lipschitz continuity of surrogate model]
\label{ass:lipschitz}
    Consider a set $\Omega \subseteq \mathbb{R}^n$ containing the origin in its interior.
    Let the surrogate dynamics~\eqref{eq:model} be locally Lipschitz continuous in the first argument on~$\Omega$ uniformly in $u \in \mathbb{U}$ and $\varepsilon \in (0, \bar{\varepsilon}]$, i.e., for every compact set $K \subseteq \Omega$, there exists $\bar{L}  = \bar{L}(K)$ such that, for every $\varepsilon \in (0,\bar{\varepsilon}]$, there is a Lipschitz constant $\LF$ with $\LF \leq \bar{L}$ satisfying
    \begin{align*}
        \|f^\varepsilon(x,u) - f^\varepsilon(y,u)\| \leq \LF\|x-y\| \quad\forall\,x,y \in K, u \in \mathbb{U}.
    \end{align*}
\end{assumption}
\noindent In Section~\ref{sec:kEDMD}, we verify Assumptions~\ref{ass:error_bound} and~\ref{ass:lipschitz} for data-driven surrogate models~\eqref{eq:model} generated by kernel EDMD using Koopman operator theory under suitable smoothness assumptions on the system under control on bounded sets~$\Omega$.

In the following MPC scheme, we use quadratic stage costs
\begin{align}\label{eq:stagecosts}
    \ell(x, u) = \|x\|_Q^2 + \|u\|_R^2 = x^\top Q x + u^\top R u
\end{align} 
with symmetric and positive-definite weighting matrices $Q \in \R^{n \times n}$ and $R\in \R^{m \times m}$.
Moreover, we require the following notion of admissibility. 
\begin{definition}[Admissible control sequences]
\label{def:admissibility}
    A control sequence $\mathbf{u} = (u(k))_{k=0}^{N-1} \subset \mathbb{U}$ of length~$N \in \mathbb{N} \cup \{+\infty\}$ is said to be \textit{admissible}. 
    The set of admissible control sequences is denoted by~$\mathcal{U}_N$. 
\end{definition}
Admissibility according to Definition~\ref{def:admissibility} does not depend on the system dynamics meaning that the definition is the same for the original dynamics~\eqref{eq:system} and the surrogate model~\eqref{eq:model}.

The proposed data-driven MPC scheme is summarized in Algorithm~\ref{alg:MPC}, where $x_{\mu_N^\varepsilon}(k)$ denotes the state of the original system governed by the dynamics~\eqref{eq:system} at time $k$ under the MPC feedback law $\mu_N^\varepsilon:\R^n \rightarrow \R^m$ computed by solving~\eqref{eq:OCP} using the data-driven surrogate dynamics~$f^\varepsilon$.
\begin{algorithm}[htb]
    \caption{Data-driven MPC}\label{alg:MPC}
    \raggedright
    \smallskip\hrule
    \smallskip
    {\it Input:} Horizon $N \in \N$, 
    surrogate~$\F$, stage cost $\ell$, input constraints~$\mathbb{U}$.
    \smallskip\hrule
    \medskip
    \textit{Initialization}: Set $k = 0$.\\[2mm]
    \noindent\textit{(1)} Measure current state $x_{\mu_N^\varepsilon}(k)$ and set $\hat{x} = x_{\mu_N^\varepsilon}(k)$.\\[1mm]
    \noindent\textit{(2)} Solve the optimal control problem~
    \begin{align}\label{eq:OCP}\tag{OCP}
    \begin{split}
        \min_{\mathbf{u} \in \mathcal{U}_N}(\hat{x}) \quad & \sum\nolimits_{i = 0}^{N-1} \ell(x_{\mathbf{u}}^\varepsilon(i), u(i)) \\
        \text{s.t.} \quad & x_{\mathbf{u}}^\varepsilon(0) = \hat{x} \text{ and, for $i \in [0:N-1]$,} \\
        & x_{\mathbf{u}}^\varepsilon(i+1) = \F(x_{\mathbf{u}}^\varepsilon(i), u(i))
    \end{split}
    \end{align}
    \hspace*{5mm} to obtain the
    optimal control sequence~$(u^\star(i))_{i = 0}^{N - 1}$. \\[1mm]
    \noindent\textit{(3)} Apply the MPC feedback law~$\mu_N^\varepsilon(\hat{x}) = u^\star({0})$ at the plant to generate the closed loop
    \[
        x_{\mu_N^\varepsilon}(k+1) = f(x_{\mu_N^\varepsilon}(k),\mu_N^\varepsilon(x_{\mu_N^\varepsilon}(k)))
    \]
    \hspace*{5mm} and shift $k = k + 1$, and go to Step~(1).
    \smallskip\hrule
\end{algorithm}

The (optimal) value function~$\V_N: \Omega \rightarrow \mathbb{R}_{\geq 0} \cup \{ \infty \}$ associated to the optimal control problem~\eqref{eq:OCP} is defined by 
\begin{align*}
    \V_N(\hat{x}) := \inf_{\mathbf{u} \in \mathcal{U}_N} \J_N(\hat{x},\mathbf{u})
\end{align*} 
with $\J_N(\hat{x},\mathbf{u}) := \sum_{k=0}^{N - 1} \ell(x^\varepsilon_{\mathbf{u}}(k),u(k))$ and $x^\varepsilon_{\mathbf{u}}(k)$ defined in Step~(2) of Algorithm~\ref{alg:MPC}. 
Analogously, we define the \textit{nominal} cost $J_N(\hat{x},\mathbf{u}) := \sum_{k=0}^{N - 1} \ell(x_{\mathbf{u}}(k),u(k))$, where $x_{\mathbf{u}}(k)$, $k\in [0:N-1]$, is obtained by propagating the original system dynamics \eqref{eq:system} starting from $x_{\mathbf{u}}(0) = \hat{x}$. The \textit{nominal} value function is defined as $V_N(\hat{x}) := \inf_{\mathbf{u} \in \mathcal{U}_N} J_N(\hat{x},\mathbf{u})$.

\section{Asymptotic stability of MPC}\label{sec:MPC_stability}
\noindent 
In this section, we prove asymptotic stability of the origin w.r.t.\ the MPC closed-loop dynamics defined in Algorithm~\ref{alg:MPC}. 
A common method to ensure closed-loop stability is the use of suitably constructed terminal conditions, see, e.g., \cite{chen1998quasi,magni2001stabilizing}. 
However, the design of terminal ingredients can be a demanding task for nonlinear systems. 
A possible alternative to achieve asymptotic stability without terminal conditions relies on cost controllability in combination with a sufficiently long prediction horizon~$N$ to ensure a relaxed Lyapunov inequality
\begin{equation}\nonumber
    V_N(f(\hat{x},\mu_N(\hat{x}))) \leq V_N(\hat{x}) - \alpha_N \ell(\hat{x},\mu_N(\hat{x}))
\end{equation}
with $\alpha_N \in (0,1]$, where $\mu_N$ denots the MPC control law using the original dynamics~$f$ in the optimization step~(2), see~\cite{Wort11,CoroGrun20} and also~\cite{KohlZeil23} for an extension to detectable stage costs. 
Then, asymptotic stability of the MPC closed loop can be concluded using relaxed dynamic programming~\cite{grune2008infinite}. 
$\alpha_N \in (0,1]$ is called suboptimality degree (performance index) since the MPC closed-loop costs on the infinite horizon are bounded by $V_\infty(\hat{x})/\alpha_N$, e.g., a factor two for $\alpha_N = 0.5$.
The computation of the suboptimality degree~$\alpha_N$ using cost controllability was originally proposed in~\cite{GrimMess05,TunaMess06} as well as further elaborated in~\cite{grune2009analysis,GrunPann10}, see also~\cite{Wort12} for a unifying comparison.
We follow the recent formulation of cost controllability proposed in~\cite[Assumption~1]{CoroGrun20}.
\begin{definition}[Cost controllability on a set~$S$] \label{def:cost:controllability}
    Consider a set $S \subseteq \mathbb{R}^n$ containing the origin in its interior. 
    The system~\eqref{eq:system} with stage cost~\eqref{eq:stagecosts} is \textit{cost controllable} on the set~$S$ if there exists a monotonically increasing and bounded sequence $(B_{N})_{N \in \mathbb{N}_0} = (B_{N}(S))_{N \in \mathbb{N}_0}$ such that for every $\hat{x} \in S$ there exists a control sequence $\mathbf{u} = \mathbf{u}(\hat{x}) \in \mathcal{U}_{\infty}$ that renders the set $S$ invariant and satisfies the growth bound\footnote{The first inequality follows from the definition of the value function and is only stated to link~\eqref{eq:growthbound:modified} to cost controllability.}
    \begin{align}\label{eq:growthbound:modified}
        V_{N}(\hat{x}) \leq J_{N}(\hat{x}, \mathbf{u}) \leq B_{N} \ell^\star(\hat{x})
    \end{align}
    for all $N \in \mathbb{N}$ with $\ell^\star(\hat{x}) := \min_{u \in \mathbb{U}}\ell(\hat{x}, u) = \|\hat{x}\|_Q^2$.
\end{definition}
\noindent Cost controllability links controllability to the performance measure~$\ell$ and, thus, ensures that the the system can be controlled to the origin sufficiently fast if only the prediction horizon~$N$ is chosen large enough.

Next, we derive asymptotic stability of the origin w.r.t.\ the MPC closed loop using the surrogate~\eqref{eq:model} in the optimization step. 
To this end, we show in Proposition~\ref{prop:cost:controllability:modified} that the growth condition~\eqref{eq:growthbound:modified} is essentially preserved if only the approximation is sufficiently accurate due to~\eqref{eq:ass:error_bound:proportional} and~\eqref{eq:ass:error_bound:uniform} and the assumed Lipschitz continuity of the surrogate model. 
Then, we conclude asymptotic stability in Theorem~\ref{thm:AS} and, finally, even show cost controllability in Corollary~\ref{cor:cost:controllability:surrogate}.
\begin{proposition}\label{prop:cost:controllability:modified}
    Suppose that system~\eqref{eq:system} is cost controllable on a set~$S$. 
    Further, let $\bar{N} \in \mathbb{N}$ be given and the Assumptions~\ref{ass:error_bound} and~\ref{ass:lipschitz} hold on a set~$\Omega$ satisfying $S \oplus \mathcal{B}_{r(\bar{N}) \eta^\varepsilon}(0) \subseteq \Omega$ with $r(\bar{N}) := \sum_{i=0}^{\bar{N}-1} \bar{L}^i$ for all $\varepsilon \in (0,\bar{\varepsilon}]$.
    Then, for the surrogate model~\eqref{eq:model}, the growth bound~\eqref{eq:growthbound:modified} is satisfied on~$S$ for all $N \in [1:\bar{N}]$ uniformly in $\varepsilon$, i.e., there exists a monotonically increasing sequence $(B_{N}^\varepsilon)_{N \in [1:\bar{N}]}$ parametrized in~$\varepsilon$ such that, for each pair $(\hat{x},N) \in S \times [1:\bar{N}]$, there exists $\mathbf{u} = \mathbf{u}(\hat{x}) \in \mathcal{U}_{N}$ satisfying
    \begin{align} \label{eq:growthbound:surrogate}
        \V_{N}(\hat{x}) \leq \J_{N}(\hat{x}, \mathbf{u}) \leq \B_{N} \cdot \ell^\star(\hat{x}). 
    \end{align}
    Moreover, $\lim_{\varepsilon \searrow 0} \B_{N} = B_{N}$ for all $N \in [1:\bar{N}]$. 
    The statement holds also upon switching the roles of $f$ and $f^\varepsilon$, i.e. if the cost controllability condition~\eqref{eq:growthbound:modified} holds for the surrogate model $f^\varepsilon$, then the growth bound~\eqref{eq:growthbound:surrogate} holds for the original system dynamics~$f$ for all $N \in [1 : \bar{N}]$.
\end{proposition}
\noindent Proposition~\ref{prop:cost:controllability:modified} shows that the growth bound~\eqref{eq:growthbound:modified} characterizing cost controllability of the system~\eqref{eq:system} is preserved for the surrogate model~\eqref{eq:model} up to a maximal horizon length~$\bar{N}$. 
The restriction to an arbitrary, but fixed maximal horizon is necessary to apply the proposed proof technique in view of the condition $S \oplus \mathcal{B}_{r(\bar{N}) \eta^\varepsilon}(0) \subseteq \Omega$ for all $\varepsilon \in (0,\bar{\varepsilon}]$. Otherwise, one has either to enlarge the set~$\Omega$, on which Assumptions~\ref{ass:error_bound} and~\ref{ass:lipschitz} hold, or reduce~$\bar \varepsilon$.

\noindent\textbf{Proof}. 
    For given $N \in \mathbb{N}$, let $\hat{x} \in S$ and $\mathbf{u} \in \mathcal{U}_{N}$ satisfy the growth bound~\eqref{eq:growthbound:modified}. 
    First, we show $x_{\mathbf{u}}^\varepsilon(k) \in \Omega$, $k \in [1:N]$. 
    Note that cost controllability yields $x_{\mathbf{u}}(k) \in S$, $k \in [1:\bar{N}]$. 
    We show $\| x_{\mathbf{u}}^\varepsilon(k) - x_{\mathbf{u}}(k) \| \leq r(k)\eta^\varepsilon = \sum_{i=0}^{k-1} L_{f^\varepsilon}^i \eta^\varepsilon$ by induction. For $k=1$, we exploit \eqref{eq:ass:error_bound:uniform} to derive
    \begin{align*}
        \| x_{\mathbf{u}}^\varepsilon(1) - x_{\mathbf{u}}(1) \| = \| f(\hat{x}, u(0)) - f^\varepsilon(\hat{x}, u(0)) \| \leq \eta^\varepsilon
    \end{align*}
    Now, let the assertion hold for $k$. Then, we have
    \begin{align*}
        & \| x_{\mathbf{u}}^\varepsilon(k+1) - x_{\mathbf{u}}(k+1) \| \\
        =& \| f^\varepsilon(x_{\mathbf{u}}^\varepsilon(k), u(k)) \mp f^\varepsilon(x_{\mathbf{u}}(k), u(k)) - f(x_{\mathbf{u}}(k), u(k)) \| \\
        \leq & \eta^\varepsilon + L_{f^\varepsilon} \|x_{\mathbf{u}}^\varepsilon(k) - x_{\mathbf{u}}(k)\| \leq 
        \sum\nolimits_{i=0}^{k} L_{f^\varepsilon}^i \eta^\varepsilon,
    \end{align*}
    i.e., we have $x_{\mathbf{u}}^\varepsilon(k) \in S \oplus \mathcal{B}_{r(k)\eta^\varepsilon}(0) \subseteq S \oplus \mathcal{B}_{r(\bar{N})\eta^\varepsilon}(0) \subseteq \Omega$ for all $k \in [0:\bar N]$, so that Assumptions~\ref{ass:error_bound} and~\ref{ass:lipschitz} hold.
    
    {Next, study the difference between $J_N(\hat{x}, \mathbf{u})$ and $J_N^\varepsilon(\hat{x}, \mathbf{u})$ using the notation $\underline{\lambda} := \operatorname{min} \{ \lambda_{\min}(Q), \lambda_{\min}(R) \}$ and $\bar{\lambda} := \lambda_{\max}(Q)$.
    We have 
    \begin{eqnarray}
        & & \underbrace{\ell(x_{\mathbf{u}}^\varepsilon(k), u(k))}_{= \| x_{\mathbf{u}}^\varepsilon(k) \mp x_{\mathbf{u}}(k) \|^2_Q + \|u(k)\|^2_R} - \ell(x_{\mathbf{u}}(k), u(k)) \label{eq:cost:controllability:proof0:modified} \\
        & \leq & \bar{\lambda} \| x_{\mathbf{u}}^\varepsilon(k) - x_{\mathbf{u}}(k) \|^2 + 2 \bar{\lambda} \| x_{\mathbf{u}}^\varepsilon(k) - x_{\mathbf{u}}(k) \| \|x_{\mathbf{u}}(k)\|. \nonumber
    \end{eqnarray}
    Then, using the shorthand notation $e_{k} := \| x_{\mathbf{u}}^\varepsilon(k) - x_{\mathbf{u}}(k) \|$ for the error, summing up this inequality over $k \in [0 : N - 1]$ yields (as $e_k = 0$ for $k = 0$)
    \begin{align*}
         \J_N(\hat{x}, \mathbf{u}) &\leq J_N(\hat{x}, \mathbf{u}) + \bar{\lambda} \Big[ \sum\nolimits_{k = 1}^{N - 1} e_k^2 + 2 e_k\|x_{\mathbf{u}}(k)\| \Big].
    \end{align*}
    Next, we derive bounds on $e_k^2$ and $e_k\|x_{\mathbf{u}}(k)\|$. To this end, the error $e_{k+1} := \| x_{\mathbf{u}}^\varepsilon(k + 1) - x_{\mathbf{u}}(k + 1) \|$ at time~$k+1$ is estimated using the triangle inequality by 
    \begin{align}
        & e_{k + 1} = \| \F(x_{\mathbf{u}}^\varepsilon(k), u(k)) \pm f(x_{\mathbf{u}}^\varepsilon(k), u(k)) - f(x_{\mathbf{u}}(k), u(k))\| \nonumber \\
        & \phantom{e_{k + 1}}  \leq c_x^\varepsilon \| x_{\mathbf{u}}^\varepsilon(k) \pm x_{\mathbf{u}}(k)\| + c_u^\varepsilon\|u(k)\| + L_f\| x_{\mathbf{u}}^\varepsilon(k) - x_{\mathbf{u}}(k) \| \nonumber \\
        & \phantom{e_{k + 1}}  \leq \bar{c}(\|x_{\mathbf{u}}(k)\| + \|u(k)\|) + (L_f + c_x^\varepsilon) e_k\label{eq:cost:controllability:proof1:modified} \\
        & \phantom{e_{k + 1}} \leq \bar{c} \sum\nolimits_{i=0}^k (L_f+c_x^\varepsilon)^{k-i} (\|x_{\mathbf{u}}(i)\| + \|u(i)\|). \label{eq:cost:controllability:proof2:modified}
    \end{align}
    with $\bar{c}:= \max\{c_x^\varepsilon, c_u^\varepsilon\}$. Let $d:=L_f + c_x^\varepsilon$.
    Using inequality~\eqref{eq:cost:controllability:proof1:modified} and the fact that $(a + b)^2 \leq 2 a^2 + 2 b^2$, we get 
    \begin{eqnarray*}
        e_k^2 & \leq & 4 \bar{c}^2 (\|x_{\mathbf{u}}(k-1)\|^2 + \|u(k-1)\|^2) + 2 d^2 e^2_{k - 1} \\
        & \leq & 4 \bar{c}^2 \underline{\lambda}^{-1} \cdot \ell(x_{\mathbf{u}}(k-1), u(k-1)) + 2 d^2 e^2_{k - 1} \\
        & \leq & 4 \bar{c}^2 \underline{\lambda}^{-1} \cdot \sum\nolimits_{i = 0}^{k - 1} (2 d^2)^{k - 1 - i} \ell(x_{\mathbf{u}}(i), u(i)) \\
        & = & 4 \bar{c}^2 \underline{\lambda}^{-1} \cdot \sum\nolimits_{j = 1}^{k} (2 d^2)^{j-1} \ell(x_{\mathbf{u}}(k-j), u(k-j))
    \end{eqnarray*}
    Analogously, leveraging inequality~\eqref{eq:cost:controllability:proof2:modified} and the fact that $2 \|a\| \|b\| \leq \|a\|^2 + \|b\|^2$ yields
    \begin{eqnarray*}
        & & e_k \|x_{\mathbf{u}}(k)\| \\[-4mm]
        & \leq & \bar{c} \sum\nolimits_{i=0}^{k-1} d^{k-1-i} \Big( \overbrace{\|x_{\mathbf{u}}(i)\| \|x_{\mathbf{u}}(k)\| + \|u(i)\| \|x_{\mathbf{u}}(k)\|}^{\leq \frac 12 ( \|x_{\mathbf{u}}(i)\|^2 + \|u(i)\|^2 + 2 \|x_{\mathbf{u}}(k)\|^2 )} \Big) \\
        & \leq & \frac{\bar{c}}{2\underline{\lambda}} \sum\nolimits_{i=0}^{k-1} d^{k-i-1} \Big( \ell(x_{\mathbf{u}}(i), u(i)) + 2 \ell^\star(x_{\mathbf{u}}(k)) \Big) \\
        & = & \frac{\bar{c}}{2\underline{\lambda}} \sum\nolimits_{j = 1}^k d^{j-1} \Big( \ell(x_{\mathbf{u}}(k-j), u(k-j)) + 2 \ell^\star(x_{\mathbf{u}}(k)) \Big)
    \end{eqnarray*}
    We begin with studying the term $\sum_{k=1}^{N-1} e_k^2$. Then, applying the derived inequality for $e_k^2$ yields
    \begin{equation}\nonumber
        \sum_{k = 1}^{N - 1} e_k^2 \leq \frac{4\bar{c}^{2}}{\underline{\lambda}} \sum_{k = 1}^{N - 1} \sum_{j = 1}^{k} (2 d^2)^{j-1} \ell(x_{\mathbf{u}}(k-j), u(k-j))        .
    \end{equation}
    Then, noting that the summation $\sum_{k = 1}^{N - 1} \sum_{j = 1}^{k}$ is equivalent to $\sum_{1 \leq j \leq k \leq N-1} = \sum_{j = 1}^{N-1} \sum_{k = j}^{N-1}$ and invoking the assumed cost controllability~\eqref{eq:growthbound:modified} leads to
    \begin{eqnarray*}
        \sum_{k = 1}^{N - 1} e_k^2 & \leq & \frac{4\bar{c}^{2}}{\underline{\lambda}} \sum_{j = 1}^{N-1} \Big[ (2 d^2)^{j-1} \underbrace{\sum_{k = j}^{N-1} \ell(x_{\mathbf{u}}(k-j), u(k-j))}_{= \sum_{r = 0}^{N-j-1} \ell(x_{\mathbf{u}}(r), u(r))} \Big] \\
        & \leq & \bar{c}^{2} \frac{4}{\underline{\lambda}} \sum\nolimits_{j = 1}^{N-1} (2 d^2)^{j-1} B_{N-j} \ell^\star(\hat{x}) =: \bar{c}^{2} c_{N} \ell^\star(\hat{x}) 
    \end{eqnarray*}
    We proceed similarly for the term $\sum_{k = 1}^{N - 1} e_k \|x_{\mathbf{u}}(k)\|$:
    \begin{eqnarray*}
        & & \sum\nolimits_{k = 1}^{N - 1} e_k \|x_{\mathbf{u}}(k)\| \\
        & \leq & \frac{\bar{c}}{2\underline{\lambda}} \sum_{k = 1}^{N - 1} \sum_{j = 1}^k d^{j-1} \Big( \ell(x_{\mathbf{u}}(k-j), u(k-j)) + 2 \ell^\star(x_{\mathbf{u}}(k)) \Big) \\
        & = & \frac{\bar{c}}{2\underline{\lambda}} \sum_{j = 1}^{N-1} d^{j-1} \bigg[ \underbrace{\sum_{r = 0}^{N-j-1} \ell(x_{\mathbf{u}}(r), u(r))}_{\leq B_{N-j} \ell^\star(\hat{x})} + 2 \underbrace{\sum_{k=j}^{N-1} \ell^\star(x_{\mathbf{u}}(j))}_{= (N-j) \ell^\star(x_{\mathbf{u}}(j))} \bigg] \\
        & \leq & \bar{c} \underbrace{\frac{1}{2\underline{\lambda}} \hspace*{-0.05cm} \bigg[ \sum_{j = 1}^{N-1} d^{j-1} B_{N-j} \hspace*{-0.05cm}+\hspace*{-0.05cm} 2 B_N \hspace*{-1.5mm}\max_{j \in [1:N-1]} \hspace*{-1mm} d^{j-1} (N\hspace*{-0.05cm}-j) \bigg]}_{=: \bar{c}_{N}} \hspace*{-0.05cm}\ell^\star(\hat{x}),
    \end{eqnarray*}
    with $d = L_f + c_x^\varepsilon$ and $\bar{c} = \max \{c_x^\varepsilon, c_u^\varepsilon\}$.
    Combining the previous estimates leads to
    \begin{align}\label{eq:B_N_epsilon}
        \J_N(\hat{x}, \mathbf{u}) \hspace*{-0.05cm}\leq\hspace*{-0.05cm} \left(B_N + \bar{c}^2 c_N + \bar{c} \bar{c}_N \right)\hspace*{-0.05cm} \ell^\star(\hat{x}) =: B_N^\varepsilon \ell^\star(\hat{x}),
    \end{align}
    where, for each $N \in [1:\bar{N}]$, $B_{N}^\varepsilon \rightarrow B_{N}$ for $\bar{c} := \max\{c_x^\varepsilon,c_u^\varepsilon\} \searrow 0$ as claimed.
    Finally, note that the proof is only based on the bounds on the difference between $x_\mathbf{u}(k)$ and $x_\mathbf{u}^\varepsilon(k)$. Hence, the same reasoning can be applied to derive a bound on $J_N(\hat{x}, \mathbf{u})$ starting from $J_N^\varepsilon(\hat{x}, \mathbf{u})$.
\hfill $\square$

Note that, differently from cost controllability as stated in Definition~\ref{def:cost:controllability}, the sequence $(\B_{N})_{N \in [1:\bar{N}]}$ derived in Proposition~\ref{prop:cost:controllability:modified} is finite. 
In particular, if $\max \{ c_x^\varepsilon, c_u^\varepsilon \} > 0$, the terms $B_N^\varepsilon$ defined in \eqref{eq:B_N_epsilon} grow unboundedly for $N \to \infty$.
However, a finite sequence is sufficient to infer asymptotic stability of the origin w.r.t.\ the MPC closed loop resulting from  Algorithm~\ref{alg:MPC} for sufficiently small~$\varepsilon$ (and, thus, sufficiently small~$c_x^\varepsilon$ and~$c_u^\varepsilon$ in Assumption~\ref{ass:error_bound}). 
The existence of a bounded sequence $(B^\varepsilon_{N})_{N \in \N}$ of infinite length satisfying the characteristic growth bound~\eqref{eq:growthbound:modified} of cost controllability for the data-driven surrogate model will be proved in Corollary~\ref{cor:cost:controllability:surrogate} leveraging the then established asymptotic stability.

The results of Proposition~\ref{prop:cost:controllability:modified} are similar to~\cite[Proposition~9]{bold2024data}, but a tighter bound on the difference between $B_N$ and $B_N^\varepsilon$ is provided and it is clarified that for a given $\varepsilon > 0$, the sequence $(B_N^\varepsilon)_{N \in [1:\bar{N}]}$ can only be derived up to some horizon $\bar{N}$. 
In addition, the relation between the sets~$S$ and~$\Omega$ has been clarified. 
Moreover, it is shown that the statement is symmetric and holds upon switching the role of the system and the surrogate model. 
This property, which was not studied in~\cite{bold2024data}, can be particularly useful in practical applications, in which the real dynamics is unknown and cost controllability can only be checked on the surrogate model.

In the following theorem, we derive our first main result. 
In essence, we show that, if the nominal MPC controller asymptotically stabilizes the origin, the data-driven MPC controller also ensures asymptotic stability of the origin w.r.t.\ the resulting closed-loop dynamics. 
\begin{theorem}\label{thm:AS}
    Let the assumptions of Proposition~\ref{prop:cost:controllability:modified} hold with an optimization horizon $N = \bar{N}$ such that $\alpha_N \in (0, 1]$ holds, where $\alpha_N$ is given by\footnote{The condition $\alpha_N \in (0,1]$ is always satisfied for sufficiently large~$N$ according to~\cite{GrunPann10} if the system~\eqref{eq:system} is cost controllable with a bounded sequence $(B_N)_{N \in \mathbb{N}_0}$.}
     \begin{align} \label{eq:alpha-N}
         \alpha_N := 1 - \frac{(B_2 - 1)(B_N - 1)\prod\nolimits_{i = 3}^N(B_i - 1)}{\prod\nolimits_{i = 2}^N B_i - (B_2 - 1)\prod\nolimits_{i = 3}^N(B_i - 1)}.
    \end{align} 
    Then, there exists $\varepsilon_0 \in (0,\bar{\varepsilon}]$ such that the MPC controller of Algorithm~\ref{alg:MPC} ensures asymptotic stability of the origin on each sublevel set $V_N^{-\varepsilon}(c) := \{ x \in \mathbb{R}^n \mid V_N^\varepsilon(x) \leq c \}$ satisfying the set inclusion $V_N^{-\varepsilon}(c) \subseteq S$ for $\varepsilon \in (0,\varepsilon_0)$.
\end{theorem}

\noindent\textbf{Proof}. 
    First, we show the relaxed Lyapunov inequality
    \begin{align}\label{eq:rel-lyap-ineq}
        \V_N(\F(x, \mu^\varepsilon_N(x))) \leq \V_N(x) - \alpha^\varepsilon_{N} \ell(x, \mu^\varepsilon_N(x))
    \end{align}
    for the surrogate model~$\F$ following the proof of \cite[Theorem~10]{bold2024data}. 
    First, we infer a lower bound on the value function by
        $V_N^\varepsilon(\hat{x}) = \inf_{\mathbf{u} \in \mathcal{U}_N} J_N^\varepsilon(\hat{x}, \mathbf{u}) \geq \inf_{u \in \mathbb{U}} \ell(\hat{x}, u) = \|\hat{x}\|_Q^2 \geq \underline{\lambda} \|\hat{x}\|^2$
    with $\underline{\lambda} := \lambda_{\min}(Q) > 0$. 
    Second, we derive an upper bound on the value function $V_N^\varepsilon$ on the set~$S$ directly from the imposed (and essentially preserved) growth bound~\eqref{eq:growthbound:surrogate}. 
    Then, the suboptimality index~$\alpha^{\varepsilon}_{N}$ of the surrogate model~\eqref{eq:model} is defined analogously to Formula~\eqref{eq:alpha-N} using the sequence $(B_i^{\varepsilon})_{i=2}^N$ instead. 
    Since $B_i^{\varepsilon} \searrow B_i$ converges monotonically for $\varepsilon \searrow 0$ in view of Proposition~\ref{prop:cost:controllability:modified}, we have that $\alpha^{\varepsilon}_N \in (0,\alpha_N)$ for all $\varepsilon \in (0,\varepsilon_0)$ for some sufficiently small $\varepsilon_0 \in (0,\bar{\varepsilon}]$.
    This ensures the relaxed Lyapunov inequality~\eqref{eq:rel-lyap-ineq} with Lyapunov function~$V_N^\varepsilon$, $\varepsilon \in (0,\varepsilon_0)$, by applying \cite[Theorem~5.4]{GrunPann10} and \cite[Theorem 1]{CoroGrun20} on $V^{-\varepsilon}_N(c)$ provided $c$ is chosen sufficiently small such that $V^{-\varepsilon}_N(c) \subseteq S$.
    Hence, we established the value function $V_N^\varepsilon(\hat{x})$ as a Lyapunov function for the closed loop of the surrogate~$\F$. 

    Next, we derive novel bounds to show that $V_N^\varepsilon$ is a Lyapunov function also for the system dynamics $f$.
    We have that
    \begin{eqnarray}
        & & V_N^\varepsilon(f(\hat{x}, \mu_N^\varepsilon(\hat{x})))  
        \pm V_N^\varepsilon(f^\varepsilon(\hat{x}, \mu_N^\varepsilon(\hat{x}))) \nonumber \\
        & \leq & V_N^\varepsilon(\hat{x}) - \alpha^\varepsilon_N \ell(\hat{x}, \mu_N^\varepsilon(\hat{x}))  + V_N^\varepsilon(f(\hat{x}, \mu_N^\varepsilon(\hat{x}))) - V_N^\varepsilon(f^\varepsilon(\hat{x}, \mu_N^\varepsilon(\hat{x}))). \label{eq:VN-bound}
    \end{eqnarray}
    Next, we exploit the proportional bounds of Assumption~\ref{ass:error_bound} to derive a bound for $V_N^\varepsilon(f(\hat{x}, \mu_N^\varepsilon(\hat{x}))) - V_N^\varepsilon(f^\varepsilon(\hat{x}, \mu_N^\varepsilon(\hat{x})))$.
    To this end, we recall that $f^\varepsilon(\hat{x}, \mu_N^\varepsilon(\hat{x})) = x_{\mathbf{u}^\star}^\varepsilon(1; \hat{x})$, where $\mathbf{u}^\star$ is the optimal solution of \eqref{eq:OCP} with initial state $\hat{x}$, and we denote the real successor state of the system by $x^+ := f(\hat{x}, \mu_N^\varepsilon(\hat{x}))$. 
    $\mathbf{u}^\sharp = (u^\sharp(i))_{i=0}^{N-1}$ represents the solution of the MPC optimization problem initialized with $\tilde{x}^+ := x_{\mathbf{u}^\star}^\varepsilon(1; \hat{x})$.\footnote{Note that $\mathbf{u}^\sharp$ is never computed in practice, but needed to define $\V_N(\F(\hat{x}, \mu^\varepsilon_N(\hat{x})))$ in the relaxed Lyapunov inequality~\eqref{eq:rel-lyap-ineq}, which we are going to leverage in the remainder of the proof.} 
    Then, optimality of MPC yields 
    \begin{align}\label{eq:VN-error}
        & \underbrace{V_N^\varepsilon(f(\hat{x}, \mu_N^\varepsilon(\hat{x})))}_{= V_N^\varepsilon(x^+) \leq J_N^\varepsilon(x^+, \mathbf{u}^\sharp)} - V_N^\varepsilon(f^\varepsilon(\hat{x}, \mu_N^\varepsilon(\hat{x}))) \\
        \leq & \sum\nolimits_{i=0}^{N-1} \Big( \ell(x_{\mathbf{u}^\sharp}^\varepsilon(i; x^+), u^\sharp(i)) - \ell(x_{\mathbf{u}^\sharp}^\varepsilon(i; \tilde{x}^+), u^\sharp(i)) \Big). \nonumber
    \end{align}
    Consider now the $i$-th term of this summation
    \begin{eqnarray}
        & & \ell(x_{ \mathbf{u}^\sharp}^\varepsilon(i; x^+), u^\sharp(i)) - \ell(x_{\mathbf{u}^\sharp}^\varepsilon(i; \tilde{x}^+), u^\sharp(i)) \nonumber\\
        & = & \|x_{\mathbf{u}^\sharp}^\varepsilon(i; x^+) \|_Q^2 - \| x_{\mathbf{u}^\sharp}^\varepsilon(i; \tilde{x}^+) \|_Q^2 \nonumber\\
        & \leq & \|Q\| \| x_{ \mathbf{u}^\sharp}^\varepsilon(i; x^+) - x_{\mathbf{u}^\sharp}^{\varepsilon}(i; \tilde{x}^+ 
        ) \| \| \underbrace{x_{\mathbf{u}^\sharp}^\varepsilon (i; x^+) + x_{\mathbf{u}^\sharp}^\varepsilon(i; \tilde{x}^+)}_{\mp x_{\mathbf{u}^\sharp}^\varepsilon(i;\tilde{x}^+)} \| \nonumber\\ 
        & \leq & 2 \|Q\| \| x_{\mathbf{u}^\sharp}^\varepsilon(i;\tilde{x}^+) \| \| x_{\mathbf{u}^\sharp}^\varepsilon(i; \tilde{x}^+) - x_{\mathbf{u}^\sharp}^\varepsilon(i; x^+) \| \nonumber \\
        & & + \|Q\| \| x_{\mathbf{u}^\sharp}^\varepsilon(i; \tilde{x}^+) - x_{ \mathbf{u}^\sharp}^\varepsilon(i; x^+) \|^2, \label{eq:diff-ell}
    \end{eqnarray}
    where we have used $\|a\|_M^2 - \|b\|_M^2 = (a+b)^\top M (a-b)$ and, then, the estimate $(a+b)^\top M (a-b) \leq \|M\| \|a - b\| \|a + b\|$.
    
    In the following we derive upper bounds for the terms $\| x_{\mathbf{u}^\sharp}^\varepsilon(i; \tilde{x}^+) - x_{\mathbf{u}^\sharp}^\varepsilon(i; x^+) \|$ and $\| x_{\mathbf{u}^\sharp}^\varepsilon(i; \tilde{x}^+) \|$. 
    First, we consider the term $\| x_{\mathbf{u}^\sharp}^\varepsilon(i; \tilde{x}^+) - x_{\mathbf{u}^\sharp}^\varepsilon(i; x^+) \|$. For $i = 0$, we leverage \eqref{eq:ass:error_bound:proportional} of Assumption~\ref{ass:error_bound} to infer
    \begin{eqnarray*}
        & & \| x_{\mathbf{u}^\sharp}^\varepsilon(0; \tilde{x}^+) - x_{\mathbf{u}^\sharp}^\varepsilon(0; x^+) \| = \| \tilde{x}^+ - x^+ \| \\
        & = & \| f^\varepsilon(\hat{x}, \mu_N^\varepsilon(\hat{x})) - f(\hat{x}, \mu_N^\varepsilon(\hat{x}))\| \leq c_x^\varepsilon \| \hat{x} \| + c_u^\varepsilon \|\mu_N^\varepsilon(\hat{x})\|.
    \end{eqnarray*}
    Using Lipschitz continuity of the surrogate model~\eqref{eq:model} given by~$\F$ with Lipschitz constant $\LF$ due to Assumption \ref{ass:lipschitz}, we obtain for $i \geq 1$
    \begin{eqnarray*}
        \| x_{\mathbf{u}^\sharp}^\varepsilon(i; \tilde{x}^+) - x_{\mathbf{u}^\sharp}^\varepsilon(i; x^+) \| & \leq & \LF^i \| x_{\mathbf{u}^\sharp}^\varepsilon(0; \tilde{x}^+) - x_{\mathbf{u}^\sharp}^\varepsilon(0; x^+) \| \\
        & \leq & \LF^i (c_x^\varepsilon \| \hat{x} \| + c_u^\varepsilon \|\mu_N^\varepsilon(\hat{x})\|).
    \end{eqnarray*}
    Second, we consider the term $\| x_{\mathbf{u}^\sharp}^\varepsilon(i; \tilde{x}^+) \|$. For the growth bound 
    of $f^\varepsilon$ derived in Proposition~\ref{prop:cost:controllability:modified} and for the relaxed Lyapunov inequality \eqref{eq:rel-lyap-ineq}, we have that
    \begin{align*}
        V_N^\varepsilon(f^\varepsilon(\hat{x}, \mu_N^{\varepsilon}(\hat{x}))) & = \sum\nolimits_{i=0}^{N-1} \ell(x_{\mathbf{u}^\sharp}^\varepsilon(i; \tilde{x}^+),u^\sharp(i)) \\
        & \leq V_N^\varepsilon(\hat{x}) \leq B^\varepsilon_N \ell^\star(\hat{x}) = B^\varepsilon_N \| \hat{x} \|_Q^2.
    \end{align*}
    Since, we have the inequality $\ell(x_{\mathbf{u}^\sharp}^\varepsilon(i; \tilde{x}^+),u^\sharp(i)) \geq \| x_{\mathbf{u}^\sharp}^\varepsilon(i; \tilde{x}^+) \|_Q^2$ for all $i \in [0:N-1]$,
    we also have $\| x_{\mathbf{u}^\sharp}^\varepsilon(i; \tilde{x}^+) \|_Q^2 \leq B^\varepsilon_N \| \hat{x} \|_Q^2$.
    By exploiting standard inequalities of weighted squared norms and taking the square root, we have $\| x_{\mathbf{u}^\sharp}^\varepsilon(i; \tilde{x}^+) \| \leq \sqrt{B^\varepsilon_N \bar{\lambda}/\underline{\lambda}} \| \hat{x} \|$.
    Then, substituting this inequalities in~\eqref{eq:diff-ell},
    we get
    \begin{eqnarray*}
         & & \ell(x_{\mathbf{u}^\sharp}^\varepsilon(i; x^+), u^\varepsilon(i)) - \ell(x_{\mathbf{u}^\varepsilon}^\varepsilon(i; \tilde{x}^+), u^\varepsilon(i)) \\
         & \leq & 2 \|Q\| \sqrt{B^\varepsilon_N \bar{\lambda} / \underline{\lambda}} \|\hat{x}\| \LF^i (c_x^\varepsilon \|\hat{x}\| + c_u^\varepsilon \|\mu_N^\varepsilon(\hat{x})\|) \\
         & & + \|Q\| \LF^{2i} (c_x^\varepsilon \| \hat{x} \| + c_u^\varepsilon \|\mu_N^\varepsilon(\hat{x})\|)^2 \\
         & \leq & 2 \|Q\| \sqrt{B^\varepsilon_N \bar{\lambda} / \underline{\lambda}}  \LF^i (c_x^\varepsilon \|\hat{x}\|^2 + \underbrace{c_u^\varepsilon \|\hat{x}\| \|\mu_N^\varepsilon(\hat{x})\|}_{\leq \frac{1}{2} c_u^\varepsilon (\|\hat{x}\|^2 + \|\mu_N^\varepsilon(\hat{x})\|^2)}) \\
         & & + 2 \|Q\| \LF^{2i} \Big( (c_x^\varepsilon)^2 \|\hat{x}\|^2 + (c_u^\varepsilon)^2 \|\mu_N^\varepsilon(\hat{x})\|^2 \Big).
    \end{eqnarray*}
    Substituting this back in \eqref{eq:VN-error} and \eqref{eq:VN-bound}, we obtain
    \begin{eqnarray*}
        & & V_N^\varepsilon(f(\hat{x}, \mu_N^{\varepsilon}(\hat{x}))) \\
        & \leq & V_N^\varepsilon(\hat{x}) - \alpha^\varepsilon \ell(\hat{x}, \mu_N^\varepsilon(\hat{x})) + C_x \| \hat{x} \|^2 + C_u \| \mu_N^\varepsilon(\hat{x}) \|^2,
    \end{eqnarray*}
    where $C_x$ and~$C_u$ are given by
    \begin{align*}
        C_x := 2 \|Q\| \sum_{i=0}^{N-1} \left( \sqrt{B^\varepsilon_N \bar{\lambda} / \underline{\lambda}} \LF^i (c_x^\varepsilon + \frac{1}{2} c_u^\varepsilon) + (\LF^i c_x^\varepsilon)^2 \right)
    \end{align*}
    and $C_u := 2  \|Q\| \cdot \sum_{i=0}^{N-1} \big( \sqrt{B^\varepsilon_N \frac{\bar{\lambda}}{\underline{\lambda}}} \LF^i \frac{1}{2} c_u^\varepsilon +  (\LF^i c_u^\varepsilon)^2 \big)$.
    $C_x$ and $C_u$ can be made arbitrarily small with sufficiently small proportionality constants $c_x^\varepsilon$ and $c_u^\varepsilon$. Then, there exists a sufficiently small $\varepsilon_0$ (potentially further tightened in comparison to our first choice resulting from the computation of the lower bound~$\alpha_N^{\varepsilon_0}$ on the suboptimality index) such that $\bar{c}_x$, $\bar{c}_u$ are sufficiently small to ensure the inequality
    \[
        - \alpha^\varepsilon_N \ell(\hat{x}, \mu_N^\varepsilon(\hat{x})) + C_x \| \hat{x} \|^2 + C_u \| \mu_N^\varepsilon(\hat{x}) \|^2 < -\bar{\alpha} \| \hat{x} \|^2
    \]
    for some $\bar{\alpha} \in (0,\alpha_N^{\varepsilon_0})$ and for all $\varepsilon \in (0,\varepsilon_0)$. This completes the proof and, thus, shows asymptotic stability of the MPC closed loop based on the surrogate model~\eqref{eq:model}.
    \hfill $\square$

The theoretical guarantees provided by Theorem~\ref{thm:AS} should be interpreted in a qualitative way, i.e., asymptotic stability is preserved provided that the approximation is sufficiently accurate and the optimization horizon~$N$ sufficiently large. However, the derived formulas for the minimum horizon and the maximum error bounds will, in general, be rather conservative in practice.

Based on Theorem~\ref{thm:AS} we can now even prove cost controllability of the surrogate model, i.e., existence of a bounded sequence $(B^\varepsilon_{k})_{k \in \mathbb{N}}$ satisfying the growth bound~\eqref{eq:growthbound:modified} on the sublevel set~$V_N^{-\varepsilon}(c)$ of Theorem~\ref{thm:AS}. 
\begin{corollary}[Cost controllability of the surrogate model]\label{cor:cost:controllability:surrogate}
    Let the assumptions of Theorem \ref{thm:AS} hold. Then there exists a monotonically increasing and bounded sequence $(B_k^\varepsilon)_{k \in \mathbb{N}}$ such that cost controllability holds for the surrogate model on the set $V_N^{-\varepsilon}(c) \subseteq S$, i.e. for all pairs $(\hat{x},N) \in V_N^{-\varepsilon}(c) \times \mathbb{N}$, there exists $\mathbf{u} \in \mathcal{U}_N$ satisfying
     \begin{align*}
        \V_{N}(\hat{x}) \leq \J_{N}(\hat{x}, \mathbf{u}) \leq B^\varepsilon_{N} \ell^\star(\hat{x}).
    \end{align*}
\end{corollary}

\noindent\textbf{Proof.} The proof resembles the proof of Proposition~6 in~\cite{MullWort17}. In the chain of inequalities~\eqref{eq:rel-lyap-ineq}, see the proof of Theorem~\ref{thm:AS}, we showed that there exists a sufficiently small~$\varepsilon$ such that $\V_N(\F(x, \mu^\varepsilon_N(x))) \leq \V_N(x) - \alpha^\varepsilon_N \ell(x, \mu^\varepsilon_N(x))$ holds for all $x \in V_N^{-\varepsilon}(c)$.
Then, denoting by $\mathbf{u}_\mathrm{MPC}$ the input sequence obtained applying the MPC control law $\mu^\varepsilon_N$, we can relate the infinite-horizon optimal cost to $\alpha^\varepsilon_N$
\begin{align*}
    V_\infty^\varepsilon(\hat{x}) &\leq J_\infty^\varepsilon(\hat{x}, \mathbf{u}_\mathrm{MPC}) \leq \sum_{i=0}^\infty \ell(x^\varepsilon(i), \mu_N^\varepsilon(x^\varepsilon(i))) \\
    & \leq (\alpha^\varepsilon_N)^{-1} \sum_{i=0}^\infty V_N^\varepsilon(x^\varepsilon(i)) - V_N^\varepsilon(f(x^\varepsilon(i), \mu_N^\varepsilon(x^\varepsilon(i))) \\
    & \leq (\alpha^\varepsilon_N)^{-1} V_N^\varepsilon(\hat{x}) \leq (\alpha^\varepsilon_N)^{-1} B_N^\varepsilon \ell^\star(\hat{x}).
\end{align*}
Since $V_N^\varepsilon(\hat{x}) \leq V_{N+1}^\varepsilon(\hat{x}) \leq V_\infty^\varepsilon(\hat{x})$ for all $N \in \mathbb{N}$, it follows that
\begin{equation*}
    V_N^\varepsilon(\hat{x}) \leq (\alpha^\varepsilon_N)^{-1} B_N^\varepsilon \ell^\star(\hat{x}).
\end{equation*}
Hence, the growth bound holds with $B_k^\varepsilon := B_N^\varepsilon / \alpha^\varepsilon_N$ for all $k > N$, which shows the assertion.
\hfill $\square$

Corollary~\ref{cor:cost:controllability:surrogate} shows that cost controllability is inherited by the data-driven surrogate model~\eqref{eq:model} provided that Assumptions~\ref{ass:error_bound} and~\ref{ass:lipschitz} hold. This result extends Proposition~\ref{prop:cost:controllability:modified}, in which the growth bound is only shown for finitely many optimization horizons~$N$, $N \in [1:\bar{N}]$.

\begin{remark}
\label{rem:state-constraints}
    We point out that the sets $\Omega$ and $S$ are not employed in the MPC algorithm meaning that we have not yet incorporated state constraints.
    Rather, we assume cost controllability of the original system~\eqref{eq:system} on a set~$S$ and, then, require Assumptions~\ref{ass:error_bound} and~\ref{ass:lipschitz} on a (sufficiently large) set~$\Omega$, such that $S \oplus r(\bar N) \eta^{\varepsilon} \subseteq \Omega$, in which we have to collect data to generate the surrogate model~\eqref{eq:model}. 
    Clearly, a larger set~$S$ (and, thus, also a larger set~$\Omega$) allows for a larger sublevel set, which corresponds to the (guaranteed) domain of attraction. 
    The incorporation of state constraints is left for future research and requires additional care to ensure recursive feasibility. 
    This might, e.g., be done based on the techniques proposed in~\cite{boccia2014stability} or using the recently introduced concept of a constraint horizon~\cite{NascWang25}, where only a finite number of predicted state are forced to obey the imposed state constraints.
    Moreover, we incorporated state constraints in a data-driven MPC with stabilizing terminal conditions in our recent work~\cite{schimperna2025stability}, where the constraints were tightened based on the Lipschitz constant~$\bar{L}$ of the surrogate model. 
    Then, combining both ideas using Assumptions~\ref{ass:error_bound} and~\ref{ass:lipschitz} may resolve problems with so-far missing state constraints. 
\end{remark} 
\noindent We emphasize the symmetry in the derived results meaning that verifying cost controllability of the surrogate model~\eqref{eq:model} for sufficient accuracy parameter~$\varepsilon$ would suffice to conclude cost controllability of the original system and, thus, also to ensure asymptotic stability of the MPC closed loop based on the surrogate model, but applied to the original system~\eqref{eq:system}.

\section{Data-driven models via kernel EDMD}\label{sec:kEDMD}

\noindent In this section, we verify Assumptions~\ref{ass:error_bound} and~\ref{ass:lipschitz}, i.e., the assumptions imposed on the surrogate model~\eqref{eq:model}, using Koopman-operator theory.
To this end, we leverage kernel extended dynamic mode decomposition (kEDMD) and its extension to control systems proposed in~\cite{BoldPhil24} as a data-driven method to approximate the Koopman operator. Moreover, at the end of the section we briefly discuss kernel approximations of the original dynamics~\eqref{eq:system} to emphasize that using the Koopman operator only serves as one of potentially many options to generate data-driven surrogate models satisfying Assumptions~\ref{ass:error_bound} and~\ref{ass:lipschitz}.

Let $\k: \R^n \times \R^n \rightarrow \R$ be a continuous, symmetric and strictly positive-definite kernel function, i.e., for every set~$\mathcal{X} = \{x_1, \dots, x_d\}\subset \R^n$ of pairwise distinct elements, the \textit{kernel matrix} $K_\mathcal{X} = (\k(x_i, x_j))_{i, j = 1}^p$ is positive definite. 
For $x \in \R^n$, the \textit{canonical features}~$\Phi_x: \mathbb{R}^n \rightarrow \mathbb{R}$ of $\k$ are defined by $\Phi_x(x') = \k(x, x')$, $x' \in \R^n$. 
By completion, the kernel~$\k$ induces an Hilbert space~$\mathbb{H}$ with inner product~$\langle \cdot,\cdot \rangle_{\mathbb{H}}$, see~\cite{Wend04}. Importantly, elements $f \in \mathbb{H}$ fulfill the reproducing property $f(x) = \langle f, \Phi_x \rangle_{\mathbb{H}}$ for all $x \in \mathbb{R}^n$ showing that point evaluation is well defined in~$\mathbb{H}$.
$\mathbb{H}$ is called \emph{reproducing kernel Hilbert space}~(RKHS). 
In this work, we use piecewise-polynomial and compactly-supported kernel functions based on the Wendland radial basis functions (RBFs)~$\Phi_{n, k}: \R^n \rightarrow \R$ with smoothness degree~$k \in \N$. 
The RKHS induced by the Wendland kernels coincides with fractional Sobolev spaces with equivalent norms \cite{Wend04}, a key property, which also holds for Matérn kernels; see \cite{FassYe11}. 
The induced (Wendland) kernel is given by $\k(x,y) = \phi_{n,k}(\|x-y\|)$ for $x, y \in \R^n$.  
In the numerical simulations of Section~\ref{sec:simulations}, $k = 1$ is used, which corresponds to $\phi_{n, 1} \in \mathcal{C}^{2}([0,\infty),\mathbb{R})$ defined by
\[
    \phi_{n, 1}(r) = \begin{cases} 
        \frac{1}{20}(1-r)^4 (4r+1) & \text{for }n \in \{ 2,3 \} \\ 
        \frac{1}{30}(1-r)^5 (5r+1) & \text{for }n \in \{4,5\}
    \end{cases}
\]
for $r < 1$ (and $0$ otherwise), see also \cite[Table~9.1]{Wend04}.

\subsection*{Koopman operator and EDMD}
Let $\Omega \subset \R^n$ be an open and bounded set with Lipschitz boundary containing the origin in its interior. 
We consider the autonomous discrete-time dynamical system given by
\begin{align}\label{eq:dynamics:F}\tag{DS}
    x^+ = F(x)
\end{align}
with map~$F: \Omega \rightarrow \R^n$.
In the following, we assume that the set~$\Omega$ is forward invariant w.r.t.\ the dynamics~\eqref{eq:dynamics:F}, i.e., 
$F(x) \in \Omega$ for all $x \in \Omega$, to streamline the presentation and refer to~\cite{kohne2024infty} for the necessary (technical) modifications if this assumption does not hold.
The associated (infinite-dimensional) linear and bounded Koopman operator~$\mathcal{K}:\mathcal{C}_b(\Omega) \rightarrow \mathcal{C}_b(\Omega)$ is defined by the identity 
\begin{align*}
    (\mathcal{K} \psi)(\hat{x}) = \psi(F(\hat{x})) \qquad \forall \ \hat{x} \in \Omega, \ \psi \in \mathcal{C}_b(\Omega),
\end{align*}
where $\mathcal{C}_b(\Omega)$ is the space of continuous and bounded functions on the set~$\Omega$.
For a set of $d \in \N$ pairwise distinct points
\begin{align}\label{eq:Xcal}
    \mathcal{X} = \{x_1, \dots, x_d\} \subset \Omega,
\end{align} 
$V_\mathcal{X} := \operatorname{span}\{\Phi_{x_1}, \dots, \Phi_{x_d}\}\subset \mathbb{H}$ is the $d$-dimensional subspace of features.
Further, $P_\mathcal{X}$ denotes the $\mathbb{H}$-orthogonal projection onto $V_\mathcal{X}$, i.e., for $h \in \mathbb{H}$, $P_\mathcal{X} h$ solves the regression problem 
    $\min_{g \in V_\mathcal{X}} \|g - h\|^2_\mathbb{H}$.
By the reproducing property, $P_\mathcal{X} h$ interpolates~$h$ at $\mathcal{X}$.
Let $\mathcal{K}|_{V_{\mathcal{X}}}$ be the restriction of the Koopman operator on the subspace $V_{\mathcal{X}}$.
As introduced in \cite{kohne2024infty}, we consider the matrix approximant $\widehat{K}$ of $\mathcal{K}|_{V_{\mathcal{X}}}$, which is given by
\begin{align}\label{eq:Koopman}
    \widehat{K} = K_\mathcal{X}^{-1} K_{F(\mathcal{X})} K_\mathcal{X}^{-1} \in \R^{d \times d}, 
\end{align}
where $K_{F(\mathcal{X})} = (\k(x_i, F(x_j)))_{i, j = 1}^d$ and $K_\mathcal{X}$ the kernel matrix corresponding to $\mathcal{X}$. In the following, we tacitly associate $\widehat{K}$ with the induced linear map from $V_\mathcal{X}$ to $V_\mathcal{X}$. The approximation $\widehat{K}$ corresponds to the widely used kernel extended dynamic mode decomposition (kEDMD; \cite{williams2014kernel}).  
For an observable function $\psi \in \mathbb{H} \subseteq \mathcal{C}_b(\Omega)$, the surrogate dynamics are given by
\begin{align*}
    \psi(F(x)) \approx \psi^+(x) := \sum\nolimits_{i = 1}^d (\widehat{K} \psi_\mathcal{X})_i \Phi_{x_i}(x), 
\end{align*}
where $\psi_\mathcal{X} = (\psi(x_1), \dots, \psi(x_d))^\top$.
A pointwise bound on the full approximation error was established in~\cite[Theorem~5.2]{kohne2024infty} and is provided in the following theorem. Therein, the approximation accuracy of kEDMD can be described using the \textit{fill distance}~$h_\mathcal{X}$ of the set $\mathcal{X}$ in $\Omega$, i.e.,
\begin{align*}
    h_\mathcal{X} := \sup_{x \in \Omega} \min_{x_i \in \mathcal{X}} \|x - x_i \|.
\end{align*}
\begin{theorem}
    \label{thm:koehne}
    Let $\mathbb{H}$ be the RKHS on~$\Omega$ generated by the Wendland kernels with smoothness degree $k \in \N$. Let ${F \in \mathcal{C}^p_b(\Omega,\mathbb{R}^n)}$ with $p = \lceil \frac{n + 1}{2} + k\rceil$ hold for system~\eqref{eq:dynamics:F}. 
    Then, there exist constants $C,h_0 > 0$ such that the bound 
    \begin{align}\label{eq:error_bound}
        \| \mathcal{K} - \widehat{K}\|_{\mathbb{H} \rightarrow \mathcal{C}_b({\Omega},\R^n)} \le Ch_\mathcal{X}^{k+1 \mathbin{/} 2}
    \end{align}
    on the full approximation error holds for all sets $\mathcal{X} :=\{ x_i \mid i \in [1:d] \} \subset \Omega$, $d \in \mathbb{N}$, of pairwise-distinct data points with fill distance~$h_{\mathcal{X}}$, $h_{\mathcal{X}} \leq h_0$.
\end{theorem}

\noindent We note that, for large data sets $\mathcal{X}$, the computation of the approximation of the Koopman operator~\eqref{eq:Koopman} can be numerically instable due to bad conditioning of the kernel matrix~$K_\mathcal{X}$. 
To alleviate this, one may include regularization by replacing ${K}_\mathcal{X}^{-1}$ with $({K}_\mathcal{X}+\lambda I)^{-1}$ with regularization parameter $\lambda > 0$. 
Error bounds in the spirit of Theorem~\ref{thm:koehne} for this regularized surrogate were proven in \cite[Theorem~2.4]{BoldPhil24}.

\subsection*{Extension of kernel EDMD to control systems.}
In this section, we recap the extension of kEDMD to systems with inputs introduced in~\cite{BoldPhil24,BoldScha25}. 
The kernel-EDMD surrogate model is defined for control-affine discrete-time systems
\begin{align}\label{eq:dynamics:control}
    x^+ = f(x, u) = g_0(x) + G(x) u
\end{align}
with locally Lipschitz-continuous maps $g_0:\Omega \rightarrow \R^n$ and $G:\Omega \rightarrow \R^{n \times m}$ with Lipschitz constants $L_{g_0}, L_G$ on~$\Omega$ and control input~$u$ restricted to a compact set $\mathbb{U} \subset \R^m$ containing the origin in its interior.
We denote by $g_i: \Omega \rightarrow \R^n$ the $i$-th column of the matrix-valued map~$G$. 
Moreover, $g_0(0) = 0$ is assumed to render the origin an equilibrium of the system~\eqref{eq:dynamics:control} for $u = 0$. 
For the relation to continuous-time systems, we refer the interested reader to \cite[Rem.~4.2]{BoldPhil24} and~\cite{StraScha25}, where it is shown that the uniform and proportional error bounds presented below are preserved for fixed~$\varepsilon$ and sufficiently small time step~$\Delta t > 0$ if only the continuous-time system is control affine. Moreover, incorporating minor modifications analogously to~\cite{StraScha26} leads to a bilinear data-driven surrogate model as shown in~\cite{StraScha25}, see also~\cite{StraWort26}.

An alternative data-driven approximation of \eqref{eq:dynamics:control} that shares the advantage of flexible state-control sampling was proposed in~\cite{bevanda2024nonparametric} using an additional kernel function to express the dependency on the control input. 
Therein, however, no error bounds were provided. 
In addition to rigorous and uniform error bounds, which are crucial for provable guarantees for data-driven MPC, the approach proposed in~\cite{BoldPhil24,SchmBold25} allows to counteract numerical ill-conditioning by proposing a two-stage approximation process as shown in the following.
\begin{assumption}[Data requirements]\label{a:data}
    For radius $\rx \geq 0$, let $\mathcal{X}$ be given by~\eqref{eq:Xcal} with $x_1 = 0$. Then, for each $i \in [1:d]$, let\footnote{$\delta_{1i}$ is the Kronecker-$\delta$, i.e., $\delta_{11} = 1$ and $\delta_{1i} = 0$ for $i \in [2:d]$.} $d_i \geq m + (1 - \delta_{1i})$ and data triplets $(x_{ij}, u_{ij}, x_{ij}^+)$, $j \in [1:d_i]$, with pairwise distinct~$u_{ij}$ be given such that ${x}_{ij} \in \mathcal{B}_{r_\mathcal{X}}(x_i) \cap \Omega$, $x_{ij}^+ = F(x_{ij}, u_{ij})$, and the matrices
    \[
        U_i := \begin{bmatrix} 
            1 & \dots & 1 \\
            u_{i1} & \dots & u_{id_i}
        \end{bmatrix} \in \R^{(m+1)\times d_i}.
    \]
    satisfy $\operatorname{rank}(U_i) = m+1$ for all $i \in [1:d]$.
\end{assumption}
\noindent Assumption~\ref{a:data} requires $D := \sum_{i=1}^{d} d_i$ data triplets. 
We call $x_i$ a virtual observation (or cluster) point to emphasize that no samples at $x_i$ are assumed. Further, $\rx$ corresponds to the cluster radius of the neighborhoods around the virtual observation points~$x_i$ from which the data are sampled. 

\textbf{Step 1: Data preparation}. Approximation of $g_0(x_i)$ and $G(x_i)$ at the virtual observation points $x_i$, $i \in [1:d]$. 
The data triplets $(x_{ij}, u_{ij}, x^+_{ij})$, $j \in [1:d_i]$, are used to compute an approximation $H_i = [\tilde{g}_0(x_i) \mid \tilde{G}(x_i)] \in \mathbb{R}^{n \times (m+1)}$ of the matrix $[g_0(x_i) \mid G(x_i)]$ for each $x_i \in \mathcal{X}$ by solving the linear regression problem $\argmin_{H_i} \big\| [x^+_{i1} \mid \ldots \mid x^+_{id_i}] - H_i U_i \big\|_F$. 
The solution may be expressed by the pseudoinverse~$U_i^\dagger$, which, for every virtual observation point~$x_i$, is well defined in view of Assumption~\ref{a:data}. 
Later, we will become interested in a uniform bound on $\max_{i \in [1:d]} \sqrt{d_i} \| U_i^\dagger \|$. 
On the one hand, this may be achieved by generating (sufficiently) exciting control values~$u_{ij}$, $j \in [1:d_i]$, following the guidelines derived in~\cite{SchmBold25} leveraging, among others, subspace-angle conditions, to ensure data efficiency ($d_i = m + (1-\delta_{1i})$). On the other hand, one may also use randomly chosen input values in combination with a slightly larger~$d_i$, see \cite[Remark~4.6]{BoldPhil24}. 

To derive \emph{proportional} error bounds for the surrogate model, we incorporate our knowledge for the data point $x_1 = 0$, that is, $g_0(x_1)=0$. 
Hence, we set $\tilde{g}_0(x_1) = 0$ while $\tilde{G}(x_1)$ is obtained by solving the reduced regression problem 
\begin{align*}
    \argmin_{G_1} \Big\| [x^+_{11} \mid \ldots \mid x^+_{1d_1}] - G_1 [u_{11} \mid \cdots \mid u_{1d_1}] \Big\|_F
\end{align*}
resembling the approach presented in \cite[Section~III.B]{UmlaPohl18}.

\textbf{Step 2: Interpolation}. The interpolation coefficients are computed analogously to the autonomous case, which leads to the following propagation step of an observable~$\psi$:
\begin{equation}\label{eq:kEDMD_model}
    \psi(f(x,u)) \approx \psi^+_\varepsilon(x) := \sum_{i = 1}^d 
    \Big[ \big( \widehat{K}_0 + \sum_{j = 1}^m u_j \widehat{K}_j \big) \psi_\mathcal{X} \Big]_i \Phi_{x_i}(x) 
\end{equation}
with $\widehat{K}_j = K_\mathcal{X}^{-1}K_{\tilde{g}_j(\mathcal{X})}K_\mathcal{X}^{-1}$, where 
\[
    (K_{\tilde g_j(\mathcal{X})})_{k,l} = \k(x_k,\tilde g_j(x_l)) = \k(x_k, (H_l)_{:,j+1})
\]
for all $k,l \in [1:d]$ and all $j \in [0:m]$, where $\tilde g_j$ denotes the $j$-th column of $\tilde{G}$.
Then, we directly construct the state-space surrogate of the control-affine system~\eqref{eq:dynamics:control}, i.e., 
\begin{align}\label{eq:dynamics:kEDMD:control}
    x^+ = \F(x, u)\, =\begin{pmatrix}
        (\psi_{1})_{\varepsilon}^{+}(x, u) \\ \vdots \\ (\psi_{n})_{\varepsilon}^{+}(x, u)
    \end{pmatrix}= \g_0(x) + \G(x) u,
\end{align}
using the approximants~\eqref{eq:kEDMD_model} of the Koopman operator with observables $\psi_\ell(x) = x_\ell$ for $\ell \in \{1, \dots, n\}$, that is, 
\begin{align*}
    g_0^\varepsilon(x)_\ell &:= \sum_{i=1}^d(\widehat K_{0} (x_\ell)_\mathcal{X})_i \Phi_{x_i}(x), \\
    (G^\varepsilon(x))_{\ell,j} &:= \sum_{i=1}^d (\widehat{K}_j (x_\ell)_\mathcal{X})_i \Phi_{x_i}(x).
\end{align*}
for $\ell \in [1:n]$ and $(\ell,j) \in [1:n]\times [1:m]$, respectively. The $i$-th column of $G^\varepsilon$ will also be denoted by $g_i^\varepsilon$.

In the following, we verify Assumptions~\ref{ass:error_bound} and~\ref{ass:lipschitz} for the proposed Koopman surrogate model~\eqref{eq:dynamics:kEDMD:control} of the control-affine dynamics~\eqref{eq:dynamics:control}. To this end, we derive the inequalities~\eqref{eq:ass:error_bound:uniform} and~\eqref{eq:ass:error_bound:proportional} and establish (uniform) Lipschitz continuity on the bounded set $\Omega$. 
The bounds depend on the fill distance~$h_\mathcal{X}$ and on the cluster radius~$\rx$, and can be made arbitrarily tight by appropriately decreasing these two parameters of the proposed data-driven approximant~\eqref{eq:dynamics:kEDMD:control}. 

The uniform error bound derived in the next theorem for the surrogate \eqref{eq:dynamics:kEDMD:control} extends Theorem~A.1 of the extended arXiv version of~\cite{BoldPhil24} to control systems by adapting the arguments used in the proof of~\cite[Theorem~4.3]{BoldPhil24}.\footnote{See \url{https://arxiv.org/abs/2412.02811}.}
\begin{theorem}[Uniform approximation error]\label{thm:control_main} 
    Let $k \geq 1$ be the smoothness degree of the Wendland kernel function. Further, let the dynamics~\eqref{eq:dynamics:control} satisfy $g_0 \in \mathcal{C}_b^{\ceil{\sigma_{n,k}}}(\Omega;\R^{n})$, $G\in \mathcal{C}_b^{\ceil{\sigma_{n,k}}}(\Omega;\R^{n\times m})$ with $\sigma_{n,k} := \frac{n+1}{2} + k$. 
    Suppose that the data satisfies Assumption~\ref{a:data} with fill distance~$h_\mathcal{X}$ and cluster radius~$r_{\mathcal{X}}$ satisfying $h_\mathcal{X} \leq h_0$ for some $h_0 > 0$ and $r_{\mathcal{X}} < h_{\mathcal{X}}/2$, respectively. 
    Then, there exist constants $C_1, C_2 > 0$ such that 
    the error bound 
    \begin{equation}\label{eq:prop error:contr}
        \|f(x,u) - \F(x,u)\|_\infty \leq C_1 h_\mathcal{X}^{k+1/2} + C_2 c \|K_\mathcal{X}^{-1}\|\rx
    \end{equation}
    holds for all $(x,u) \in \Omega \times \mathbb{U}$, where $c$ is defined by
    \begin{align*}
        c := \max_{i \in [1:d]} \sqrt{d_i} \| U_i^\dagger\| \cdot \Phi_{n,k}^{1/2}(0) \max_{v \in \mathds{1}} \sqrt{v^\top {K}_\mathcal{X}^{-1}v}
    \end{align*}
    and $\mathds{1} = \{v \in \R^d \mid v_i \in \{\pm 1\}\}$.
\end{theorem}
\noindent\textbf{Proof.}
    In a first step, analogously to \cite[Proof of Theorem 4.2]{BoldPhil24}, we obtain 
    \begin{align*}
        &|H_{pq}(x_i) - (H_i)_{pq}| \leq \sqrt{2 d_i}(L_{g_0} + L_G \overline{u}) \| U_i^\dagger \| \cdot \rx 
    \end{align*}
    for $i \in [1:d]$, where 
    $H:\R^n \to \R^{n \times (m+1)}$ is defined by $H(x)=[g_0(x)\,|\,G(x)]$ and $\overline{u} := \max \{ \|u\|_\infty : u \in \mathbb{U} \}$.
    Due to uniform continuity of the kernel~$\k$ on $\overline{\Omega}\times \overline{\Omega}$, we get the estimate
    \begin{eqnarray}\label{eq:est:K}
        & & \|K_{\tilde{g}_j(\mathcal{X})} - K_{{g}_j(\mathcal{X})}\| \nonumber \\
        & \leq & {c_\k} (L_{g_0} + L_G \overline{u}) \cdot \max \{ \sqrt{2 d_i} \| U_i^\dagger\| \mid i \in [1:d] \} \rx 
    \end{eqnarray}
    for a constant ${c_\k}\geq 0$ depending on the continuity modulus of the kernel~$\k$. 
    Consequently, \eqref{eq:est:K} implies
    \begin{align*}
        \| \widehat{K}_j - K_\mathcal{X}^{-1}K_{g_j(\mathcal{X})}K_\mathcal{X}^{-1}\| &=  \| K_\mathcal{X}^{-1}(K_{\tilde{g}_j(\mathcal{X})} - K_{g_j(\mathcal{X})})K_\mathcal{X}^{-1}\|\\
        & \leq \|K_\mathcal{X}^{-1}\|^2 \|K_{\tilde{g}_j(\mathcal{X})} - K_{{g}_j(\mathcal{X})}\|\\
        & \leq \tilde{c} \|K_\mathcal{X}^{-1}\|^2 \max_{i \in [1:d]} \sqrt{d_i} \| U_i^\dagger\| \cdot \rx
    \end{align*}
    with $\tilde{c} = \sqrt{2} c_\k \cdot (L_{g_0} + L_G \overline{u})$.
    
    Next, we show inequality~\eqref{eq:prop error:contr} with $h_{\mathcal{X}}^{k-1/2} \operatorname{dist}(x, \mathcal{X})$ instead of $h_{\mathcal{X}}^{k + 1/2}$, which will turn out to be beneficial in the proof of Proposition~\ref{prop:error_bound:proportional}.
    Using \cite[Theorem 11.17]{Wend04} with multiindex~$|\alpha| = 1$ together with \cite[Theorem 3.7]{BoldPhil24}, we get the proportional bound
    \begin{eqnarray*}
        & & |(K_\mathcal{X}^{-1} K_{g_j(\mathcal{X})} K_\mathcal{X}^{-1}\varphi_\mathcal{X})^\top\textbf{k}_\mathcal{X}(x) - \varphi \circ g_j(x) | \\
        & \leq & \hat{c}h_\mathcal{X}^{k-1/2}\|\varphi\|_{\mathbb{H}} \ {\operatorname{dist}(x,\mathcal{X})}
    \end{eqnarray*}
    with $\hat{c}\geq 0$ for all $\varphi \in \mathbb{H}$ 
    and $\textbf{k}_\mathcal{X} = (\Phi_{x_1}, \dots, \Phi_{x_d})^\top$. Thus, by the triangle inequality, we obtain the following bound on the full error:
    \begin{eqnarray*}
        && |(\widehat{K}_j\varphi_\mathcal{X})^\top\textbf{k}_\mathcal{X}(x)  - \varphi \circ g_j(x) | \\
        & \leq & \Big( \hat{c}h_\mathcal{X}^{k-1/2}\operatorname{dist}(x,\mathcal{X}) + \tilde{c} c \|K_\mathcal{X}^{-1}\| \rx \Big) \|\varphi\|_{\mathbb{H}}.
    \end{eqnarray*}
    Here, the terms~$\Phi_{n,k}^{1/2}(0)$ and $\max_{v \in \mathds{1}} \sqrt{v^\top {K}_\mathcal{X}^{-1}v}$ in $c$ originate from the bound on the norm of $\varphi_\mathcal{X}$ and $\mathbf{k}_\mathcal{X}(x)$, respectively, similar to \cite[Proof of Theorem 4.2]{BoldPhil24}.
    Last, choosing the coordinate functions as observables, we may, for the remainder of the proof, assume that $\varphi$ is linear since we choose $\varphi(x) = x_i$ with finite RKHS norm $\|\varphi\|_{\mathbb H}$ due to the Sobolev characterization of the Wendland RKHS. Thus,
    \begin{align*}
        \varphi(x^+) = \varphi(g_0(x) + G(x)u) = \varphi(g_0(x)) + \sum_{j=1}^m\varphi(g_j(x))u_j
    \end{align*}
    and hence, we get \eqref{eq:prop error:contr} with $C_1 = \hat{c}\max_{i\in [1:n]}\|x_i\|_\mathbb{H}$ and $C_2 = \tilde{c} \max_{i\in [1:n]}\|x_i\|_\mathbb{H}$, where $x_i$, $i\in [1:n]$ denotes the $i$-th coordinate map.
    \hfill $\square$
    
As can be seen in Theorem~\ref{thm:control_main}, the approximation quality mainly depends on the fill distance~$h_\mathcal{X}$ of the grid~$\mathcal{X}$ of cluster points and the cluster radius $\rx \geq 0$. 
The choice of the virtual observation points~$\mathcal{X}$ not only influences the first term of the error bound via the fill distance, but also $\|K_\mathcal{X}^{-1}\|$ and~$c$ in the second term. 
Depending on the choice of~$\mathcal{X}$, the cluster radius~$\rx$ has to be chosen sufficiently small to compensate for $C_2c\|K_\mathcal{X}^{-1}\|$.
Hence, by choosing the cluster radius~$\rx$ such that the two terms are \textit{balanced}, we may ensure the uniform upper bound $\eta^\varepsilon := 2C_1 h_0^{k+1/2}$ such that~\eqref{eq:ass:error_bound:uniform} holds.

In the following lemma we verify Assumption~\ref{ass:lipschitz} for the kEDMD surrogate model, i.e., we show that the surrogate model is Lipschitz in the first argument, uniformly in $u$ and~$\varepsilon$.
\begin{lemma}[Lipschitz continuity of surrogate~$\F$]\label{lem:uniformcont} 
    Let the assumptions of Theorem~\ref{thm:control_main} hold.
    Then, there exist constants $h_0, \bar{L}$ such that, if the fill distance~$h_\mathcal{X}$ and the cluster radius~$r_{\mathcal{X}}$ satisfy $h_\mathcal{X} \leq h_0$ and $r_{\mathcal{X}} \leq h_{\mathcal{X}}/2$, the surrogate~$\F$ defined by~\eqref{eq:dynamics:kEDMD:control} is Lipschitz continuous in~$x$ uniformly in $u$ with Lipschitz constant~$\LF \leq \bar{L}$, i.e.
    \begin{align}\label{eq:kEDMD:Lipschitz}
        \| \F(x, u) - \F(y, u)\| \le \LF \|x - y\|
    \end{align}
    holds for all $x, y \in \Omega$ and $u \in \mathbb{U}$. 
\end{lemma}
\noindent\textbf{Proof.} 
    Let $x,y \in \Omega$ and $u \in \mathbb{U}$ be given. Then, it holds by using Taylor's theorem 
    \begin{align*}
        \F(x, u) & = \F(y,u) + (x - y)^\top J_{\F, x}(\xi,u)
    \end{align*}
    for some $\xi \in \{\tilde{\xi} \in \R^n \mid \tilde{\xi} = (t - 1)x + t y, \ t \in [0, 1]\}$.
    Hence, we may estimate
    \begin{eqnarray*}
        & & \|\F(x, u) - \F(y, u)\| \\
        & \leq & \|x - y\| \|J_{\F, x}(\xi, u) \pm J_{f, x}(\xi, u)\| \\
        & \leq & \|x - y\|(\|J_{\F, x}(\xi, u) - J_{f, x}(\xi, u)\|  + \|J_{f, x}(\xi, u)\|).
    \end{eqnarray*}
    In the following, let $\hat{g}_i$ for $i \in [0:m]$ be the approximation of $g_i$ obtained with $\rx = 0$, i.e. without errors in the first step of the identification algorithm.
    To compute a bound for $\|J_{\F, x}(\xi, u) - J_{f, x}(\xi, u)\|$, we first leverage the control-affine structure and, then, apply {\cite[Theorem 11.17]{Wend04} with multiindex~$|\alpha| = 1$} to obtain
    \begin{align*}
        & \|J_{\F, x}(z, u) - J_{f, x}(z, u)\| \\
        \le & \|J_{\g_0}(z) - J_{g_0}(z) \pm J_{\hat{g}_0}(z)\|  + \sum\limits_{i = 1}^m \| (J_{\g_i}(z) - J_{g_i}(z) \pm J_{\hat{g}_i}(z)) u_i\| \\
        \leq & \Big[ C_{J_{g_0}} + \max_{u \in \mathbb{U}}\ \sum_{i=1}^m C_{J_{g_i}} |u_i| \Big] h_\mathcal{X}^{k-1/2} \\
        & + \|J_{\hat{g}_0}(z) - J_{\g_0}(z) \| + \sum\limits_{i = 1}^m\|(J_{\hat{g}_i}(z) - J_{\g_i}(z)) u_i \|,
    \end{align*}
    where $C_{J_{g_i}}$, $i \in [0:m]$, are constants independent of $h_\mathcal{X}$ and $\rx$. 
    Analogously to the proof of \cite[Theorem 4.3]{BoldPhil24}, letting $u_0 = 1$, we can estimate for $i \in [0:m]$ and $j \in [1:n]$ 
    \begin{eqnarray}
        & & \| [ \nabla \hat{g}_{ij}(z) - \nabla \g_{ij}(z)] u_i \| \nonumber\\
        & \leq & \max_{u \in \mathbb{U}}|u_i| L_{J_{g_i}} 
     \sqrt{2\max_{k\in [1:d]}d_k} \left( \max_{\ell \in [1:d]} \|U^\dagger_\ell\|\right) \nonumber\\
        & & \cdot \|K_\mathcal{X}^{-1}\textbf{k}_\mathcal{X}(z)\|\|K_\mathcal{X}^{-1}\| \rx \nonumber\\
        & \leq & \max_{u \in \mathbb{U}}|u_i| L_{J_{g_i}} \sqrt{2\max_{k\in [1:d]}d_k}  \left( \max_{\ell \in [1:d]} \|U^\dagger_\ell\|\right)\Phi_{n, k}(0) \nonumber \\
        & & \cdot \max_{v \in \mathds{1}
        }v^\top K_\mathcal{X}^{-1}v \|K_\mathcal{X}^{-1}\| \rx =: c_{\mathcal{X}i}, \label{eq:cxi}
    \end{eqnarray}
    where $\textbf{k}_\mathcal{X}(z) = (\k(x_1, z), \dots, \k(x_d, z))^\top$ and 
    $L_{J_{g_i}}$ are the Lipschitz constants of the matrix-valued functions~$J_{g_i}$. 
    Since the function~$f$ is continuously differentiable w.r.t.\ its first argument on the compact set~$\Omega$, we have $\| J_{f, x}(\xi, u)\| \leq \widetilde{C} \in (0,\infty)$. 
    Overall, we have verified inequality~\eqref{eq:kEDMD:Lipschitz} with Lipschitz constant~$\LF$ given by 
    \[
         \widetilde{C} + C_{J_{g_0}} h_\mathcal{X}^{k-1/2} + m\! \max\limits_{i \in [1:m], u \in \mathbb{U}} C_{J_{g_i}} |u_i| h_\mathcal{X}^{k-1/2} + n \sum\nolimits_{i=0}^m c_{\mathcal{X}i}.
    \]
    Hence, suitably adjusting the cluster radius~$r_\mathcal{X}$ depending on the virtual observation points and, thus, $\| K_\mathcal{X}^{-1} \|$ shows the statement. 
\hfill $\square$

The continuity in Lemma~\ref{lem:uniformcont} is uniform in the approximation bound~$\varepsilon>0$. This is in contrast to finite-element dictionaries, see~\cite{schaller2023towards}, where the derivative of the ansatz functions increases for decreasing mesh size.

The uniform error estimate in Theorem~\ref{thm:control_main} still contains a constant offset in the upper bound, which results from Step~1 of the approximation scheme. 
Given an exact approximation of $g_0$ and $G$ at the grid points in $\mathcal{X}$, i.e., $\rx = 0$, 
and inspecting the last part of the proof of Theorem~\ref{thm:control_main}, one may straightforwardly deduce an error bound~\eqref{eq:prop error:contr} that is proportional in the distance to the grid points. More precisely, as $x_1=0$ due to Assumption~\ref{a:data}, the first term in the approximation error~\eqref{eq:prop error:contr} may be replaced by
\begin{align*}
    {C_1} h_\mathcal{X}^{k-1/2} \operatorname{dist}(x, \mathcal{X}) \leq {C_1} h_\mathcal{X}^{k-1/2} \|x\|.
\end{align*}
In the following, we show how to ensure a proportional error bound even for $\rx > 0$. 
To this end, we leverage the fact that we used the knowledge of the drift dynamics in the origin, i.e., we enforced $g_0^\varepsilon(0) = 0$ in the first step of the algorithm so that the model exactly describes the system in the origin. 
\begin{proposition}[Proportional error bound]\label{prop:error_bound:proportional}
    Let the assumptions of Theorem~\ref{thm:control_main} hold. 
    Then, there exist constants $h_0, r_0$ such that, if the fill distance~$h_\mathcal{X}$ and the cluster radius~$r_\mathcal{X}$ satisfy $h_\mathcal{X} \leq h_0$ and $r_\mathcal{X} \leq r_0$, we have the proportional error bound
    \begin{equation*}
        \|f(x,u) - \F(x,u)\| \leq c_x^\varepsilon \|x\| + c_u^\varepsilon \|u\|,
    \end{equation*}
    for all $x \in \Omega$, $u \in \mathbb{U}$,
    where $c_u^\varepsilon := C_1 h_\mathcal{X}^{k-1/2} \operatorname{dist}(x, \mathcal{X}) + C_2 c \|K_\mathcal{X}^{-1}\|\rx$, and 
    \begin{equation*} 
        c_x^\varepsilon := C_{J_{g_0}} h_\mathcal{X}^{k-1/2} + m \max\limits_{i \in [1:m], u \in \mathbb{U}} C_{J_{g_i}} |u_i| h_\mathcal{X}^{k-1/2} + n \sum_{i=0}^m c_{\mathcal{X}i},
    \end{equation*}
    where $c_{\mathcal{X}i}$ are defined in \eqref{eq:cxi} and depend linearly on the cluster radius $\rx$.
\end{proposition}
\noindent\textbf{Proof.} 
    To prove the statement, we first show that the difference $\| J_{\F}(x, u) - J_{f}(x, u) \|$ between the differentials of $f$ and $\F$ is bounded.
    From the proof of Lemma \ref{lem:uniformcont} we obtain the bound 
    \begin{align*}
    &\quad \| J_{\F,x}(x, u) - J_{f,x}(x, u) \| \\ & \leq C_{J_{g_0}} h_\mathcal{X}^{k-1/2} + m \max\limits_{i \in [1:m], u \in \mathbb{U}} C_{J_{g_i}} |u_i| h_\mathcal{X}^{k-1/2} + n \sum_{i=0}^m c_{\mathcal{X}i},
    \end{align*}
    i.e., $c_x^\varepsilon$. 
    To compute a bound for $\|J_{\F,u}(x, u) - J_{f,u}(x, u)\|$, we exploit the error bounds of Theorem~\ref{thm:control_main}, to obtain
    \begin{eqnarray*}
        & & \|J_{\F,u}(x, u) - J_{f,u}(x, u)\| = \|G^\varepsilon(x) - G(x)\| \\
        & \leq & C_1 h_\mathcal{X}^{k-1/2} \operatorname{dist}(x, \mathcal{X}) + C_2 c \|K_\mathcal{X}^{-1}\|\rx =: c_u^\varepsilon.
    \end{eqnarray*}
    
    Now, let $e(x, u) := \F(x,u) - f(x,u)$. In view of the modified learning algorithm, there is no approximation error in the origin, i.e. $e(0,0) = 0$. We now study $e(x,u)$ using Taylor's theorem. There exists $t \in [0,1]$ such that, denoting $z = tx$ and $v = tu$,
    \begin{eqnarray*}
        \|e(x,u)\| & = & \|\underbrace{e(0,0)}_{=0} + J_{e, x}(z,v) x + J_{e, u}(z, v) u\| \\
        & \leq & \|J_{\F,x}(z, v) - J_{f,x}(z, v)\| \|x\| + \|J_{\F,u}(z, v) - J_{f,u}(z, v)\| \|u\| \\
        & \leq & c_x^\varepsilon \|x\| + c_u^\varepsilon \|u\|,
    \end{eqnarray*}
    which proves the statement.
\hfill $\square$

\begin{remark}
    For the kEDMD surrogate, the coefficients in the proportional error bounds $c_x^\varepsilon$ and $c_u^\varepsilon$ are functions of the fill distance $h_\mathcal{X}$ and of the cluster radius $\rx$, and can be rewritten as
    \begin{align*}
    \begin{array}{rll}
        c_x^\varepsilon &= c_x^h h_\mathcal{X}^{k-1/2} &+ c_x^r \rx \\
        c_u^\varepsilon &= c_u^h h_\mathcal{X}^{k-1/2} \operatorname{dist}(x, \mathcal{X}) &+ c_u^r \rx
        \end{array}
    \end{align*}
    for suitable $c_x^h, c_x^r, c_u^h, c_u^r$. The terms $c_x^r$ and $c_u^r$ have a dependence on $\|K_{\mathcal{X}}^{-1}\|$, which is related to the location and number~$d$ of virtual observation points. 
    For increasing~$d$, the minimal eigenvalue of the kernel matrix~$K_\mathcal{X}$ may decrease, and consequently $\|K_{\mathcal{X}}^{-1}\|$ may increase. 
    Hence, to satisfy given bounds for $c_x^\varepsilon$ and $c_u^\varepsilon$, one should first select
    the location of the virtual observation points (that is related to the fill distance $h_\mathcal{X}$) to obtain sufficiently small terms $c_x^h h_\mathcal{X}^{k-1/2}$ and $c_u^h h_\mathcal{X}^{k-1/2} \operatorname{dist(x, \mathcal{X})}$. Given the location of the virtual observation points, $\|K_{\mathcal{X}}^{-1}\|$ is fixed. Then, in a second step, it is possible to select the cluster radius $\rx$ sufficiently small to satisfy the required bounds for the terms $c_x^r \rx$ and $c_u^r \rx$. 
\end{remark}

\noindent The proportional error bound derived in Proposition~\ref{prop:error_bound:proportional} may be more conservative than the bounds of Theorem~\ref{thm:control_main} further away from the origin, but converges to zero for $(x,u) \rightarrow (0,0)$, which is a key property to rigorously show asymptotic stability of the MPC, as shown in Theorem~\ref{thm:AS} and summarized in the following theorem.
\begin{theorem}
    Let Assumption~\ref{a:data} on the data hold and consider the data-driven MPC scheme summarized in Algorithm~\ref{alg:MPC} using the kEDMD surrogate model~\eqref{eq:dynamics:kEDMD:control} of the control-affine dynamics~\eqref{eq:dynamics:control} in the optimization step of Algorithm~\ref{alg:MPC}. Let $N \in \mathbb{N}$ be such that $\alpha_N$ defined in \eqref{eq:alpha-N} is in $(0, 1]$ and assume that the system is cost controllable on a set $S \subseteq \Omega \ominus \mathcal{B}_{r(N)\eta^\varepsilon}$, see Definition~\ref{def:cost:controllability}, with $\mathbf{u} \in \mathcal{U}_N$. Then, if the fill distance~$h_\mathcal{X}$ and the cluster radius~$\rx$ are sufficiently small, Algorithm~\ref{alg:MPC} asymptotically stabilizes the system~\eqref{eq:dynamics:control}.
\end{theorem}
\noindent\textbf{Proof.} 
    We have verified \eqref{eq:ass:error_bound:proportional} and \eqref{eq:ass:error_bound:uniform} of Assumption~\ref{ass:error_bound}, and Assumption~\ref{ass:lipschitz} for the kEDMD surrogate of the control-affine dynamics~\eqref{eq:dynamics:control} in Proposition~\ref{prop:error_bound:proportional}, Theorem~\ref{thm:control_main}, and Lemma~\ref{lem:uniformcont}, respectively. Hence, all assumptions of Theorem~\ref{thm:AS} hold in view of the imposed cost controllability and the assertion can be inferred.
\hfill $\square$

\subsection*{Relation to kernel approximation of control systems.}
The kEDMD approximation defined in \eqref{eq:Koopman} consists of three components: First, the observable is projected onto the dictionary $V_\mathcal{X}$ (where the $\mathbb{H}$-orthogonal projection is realized by $K_\mathcal{X}^{-1}$). 
Then, the linear dynamics~$K_{F(\mathcal{X})}$ encoding the knowledge of the propagated features, is applied. 
As the result is, in general, not contained in the finite-dimensional linear space~$V_\mathcal{X}$ (due to concatenation of the features with the nonlinear flow), the result is again projected. 
While the dynamics~\eqref{eq:Koopman} is formulated for arbitrary observables, we simply applied the resulting operator to the coordinate functions to deduce the surrogate model~\eqref{eq:dynamics:kEDMD:control}.

An alternative approximation scheme can be derived by \emph{directly} approximating the flow map~$F$ of the dynamical system~\eqref{eq:dynamics:F} via kernel regression, i.e., 
\begin{equation}\label{eq:dynamics:kernel}
    \widetilde F_i(x) = \sum_{j=1}^d (K_\mathcal{X}^{-1} (F_i)_\mathcal{X})_j \Phi_{x_j}(x) = (F_i)_\mathcal{X}^\top K_\mathcal{X}^{-1} \hspace*{-1mm}\begin{bmatrix} 
        \Phi_{x_1}(x) \\
        \vdots \\ 
        \Phi_{x_d}(x) 
    \end{bmatrix}
\end{equation}
for each $i \in [1:n]$. 
This corresponds to direct interpolation of the flow map and hence involves only one projection. Moreover, the main difference to the approximation \eqref{eq:Koopman} are the data requirements: For the approximation \eqref{eq:Koopman} applied to the coordinate functions, we require data through the lens of the features, i.e., $\k(x_i,F(x_j))$ for $(i,j)\in [1:d]^2$, whereas for the approximation~\eqref{eq:dynamics:kernel}, we require samples of the flow map, i.e., $F(x_j)$ for all $j\in [1:d]$. 
The error bound for the surrogate dynamics~\eqref{eq:dynamics:kernel}, however, is very closely related to the one inferred from~\eqref{eq:Koopman} and of similar structure, see \cite[Theorem 3.4]{kohne2024infty}, where the first estimate in this reference corresponds to~\eqref{eq:dynamics:kernel} (involving only one projection) and the second estimate to~\eqref{eq:Koopman} (involving two projections). As both projection errors are proportional to the fill distance~$h_\mathcal{X}$, the convergence rate is same and the results above can be derived analogously. 
Hence, Assumption~\ref{ass:error_bound} can also be verified for the surrogate model~\eqref{eq:dynamics:kernel}.

\section{Numerical simulations}\label{sec:simulations}

\noindent In this section, we illustrate the asymptotic stability of the origin w.r.t.\ the MPC closed loop using the kernel-EDMD surrogate in the prediction and optimization step by conducting numerical simulations.
In the following, we denote by \textit{physics-informed kEDMD} (\textbf{PI-kEDMD}) the models obtained with the modified Step~1 for the encoding of the equilibrium, and with \textbf{kEDMD} the models obtained according to \cite{BoldPhil24} without modifying Step~1.

\subsection*{Van der Pol}
We consider the Euler discretization of the nonlinear van-der-Pol oscillator. Due to using the Euler discretization, this directly yields the control-affine system
\begin{align}\label{eq:van der Pol}
    x^+ = x + \Delta t \binom{x_2}{\nu(1 - x_1)^2x_2 - x_1 + u} 
\end{align}
with parameters $\Delta t = 0.05$, $\nu = 0.1$, control constraints $\mathbb{U} = [-2, 2]$ and domain $\Omega = [-2, 2]^2$. 
The kEDMD surrogate~\eqref{eq:dynamics:kEDMD:control} is constructed using Wendland kernels with smoothness degree~$k = 1$.
The cluster points are chosen according to a Padua grid. The set of Padua points of order $p$ in the domain $[-1, 1] \times [-1, 1]$ is defined as
\begin{equation*}
    \{ (\mu_j, \eta_k), 0 \leq j \leq p, 1 \leq k \leq \big\lfloor \frac{p}{2} \big\rfloor + 1 + \delta_j \}
\end{equation*}
where $\delta_j = 1$ if $p$ and $j$ are both odd and $\delta_j = 0$ otherwise. Then, the pair $(\mu_j,\eta_k)$ is defined by
\begin{align*}
    \mu_j := \cos\Big(\frac{j \pi}{p}\Big), \quad \eta_k := \begin{cases}
        \cos\Big(\frac{(2k - 2) \pi}{p + 1}\Big) \quad j \text{ odd}, \\
        \cos\Big(\frac{(2k - 1) \pi}{p + 1}\Big) \quad j \text{ even}.
    \end{cases}
\end{align*}
Moreover, we add a cluster point in the origin to the Padua grid, according to Assumption \ref{a:data}. The Padua grid is properly scaled to cover the domain $\Omega$.
Padua grids of order $p \in \{ 25, 50 \}$ are considered, resulting in $d \in \{ 352, 1327 \}$ state-space sample points.
In comparison to a uniform grid, the use of a Padua grid allows to have a more uniform error through the domain, and lower condition numbers in the model matrices. 
In fact, Padua points minimize, analogously to Chebychev nodes in one dimension, the Lebesgue constant and, thus, the condition number and are, thus, preferable. The use of a uniform grid with a large number of data points may lead to numerical problems in the solution of the optimal control problem.
For each cluster point~$x_i \in \mathcal{X}$, $d_i = 25$ random control values $u_{ij} \in \mathbb{U}$ are chosen, yielding the data points~$(x_{ij}, u_{ij}, x_{ij}^+)$. 
Herein, $x_{ij}$ and $u_{ij}$ are drawn according to our data requirements with $\rx = \nicefrac{\sqrt{2}}{d}$ to specify the neighborhood of~$x_i$, see Section~\ref{sec:kEDMD}. For the derivation of the PI-kEDMD model Step 1 of the algorithm is modified for the cluster point $x_1 = 0$ by setting $\tilde{g}_0(x_1) = 0$, while in the kEDMD model Step 1 is performed in the same way for all cluster points.
The MPC cost function is chosen as in~\eqref{eq:stagecosts} with matrices $Q = I_2$ and $R = 10^{-4}$. 
The following closed-loop simulations emanate from the initial state $x^0 = (0.5, 0.5)^\top$, and consider a prediction horizon $N \in \{ 10, 30 \}$ in MPC. The simulations are implemented in Python, and Casadi is used for the solution of the MPC optimization problem.
The comparison of the closed loop error with different number of clusters and different prediction horizons is reported in Fig.~\ref{fig:vdp_comparison}.
\begin{figure}[htb]
    \centering
    \includegraphics[width=0.4\linewidth]{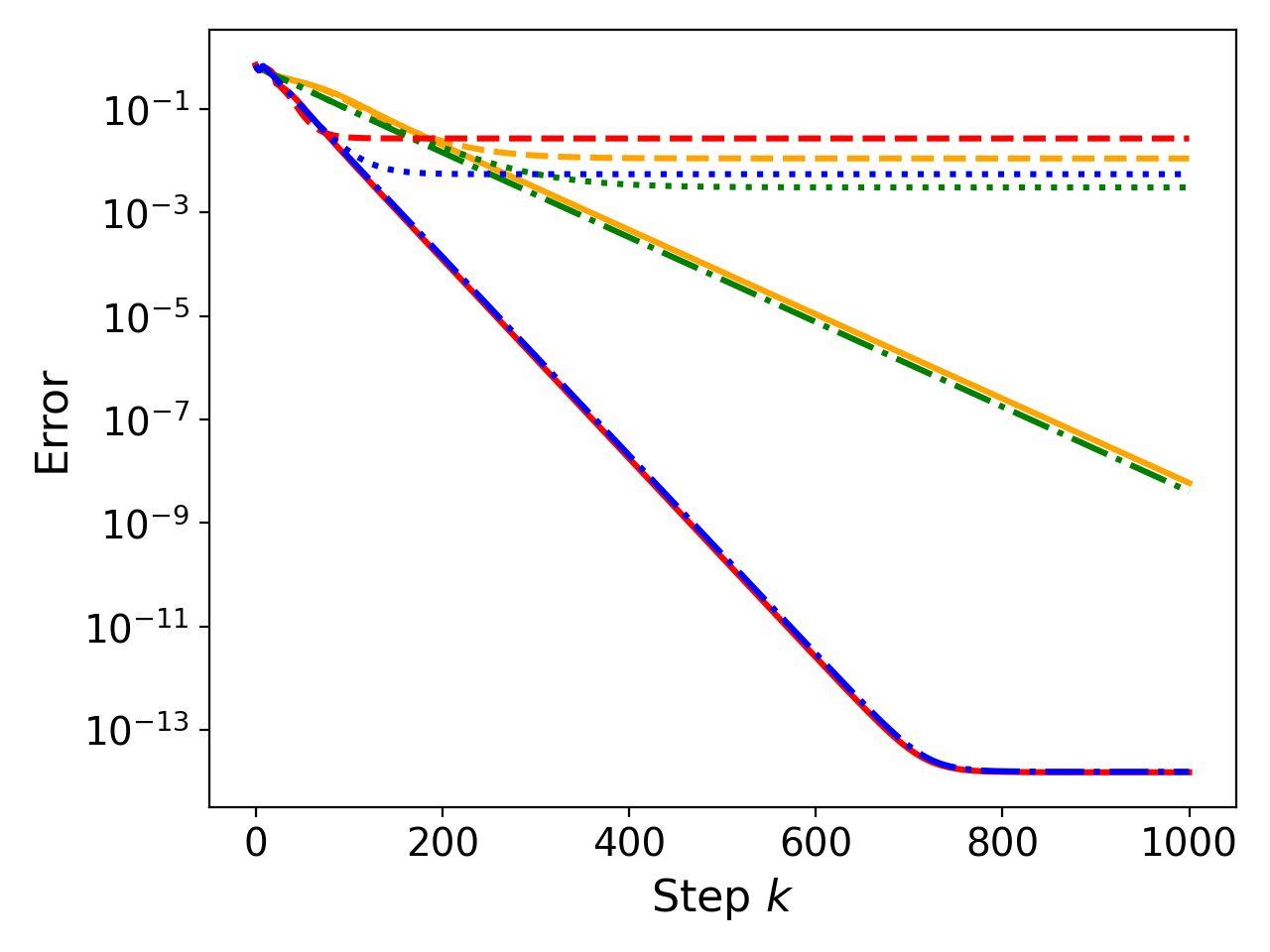}
    \includegraphics[width=0.4\linewidth]{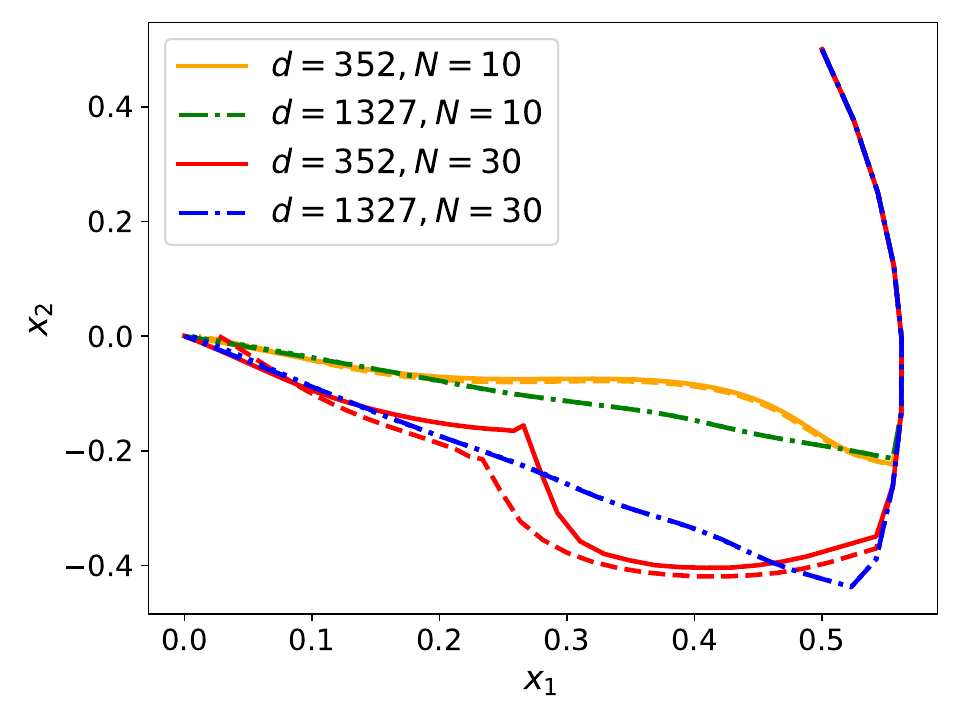}
    \caption{Van der Pol: Error $\|x(k)\|$ for the kEDMD-based MPC closed loop (left) and trajectories with different horizons and number of clusters (right). Solid and dash-dotted lines are used for \textbf{PI-kEDMD} and dashed and dotted lines for \textbf{kEDMD}.}
    \label{fig:vdp_comparison}
\end{figure}

In the simulations with the kEDMD model, we observe practical asymptotic stability only: the error stagnates at a positive constant value. 
The tracking accuracy increases with the number of cluster points. This is expected, as a large number of cluster points is related to lower modeling errors.
Instead, the PI-kEDMD models allow to reach a lower error, with a decrease up to the solver tolerance at around $10^{-14}$, in accordance to the shown asymptotic stability property of the MPC closed loop. 
We also notice that a fast decrease to the origin is obtained with the prolonged prediction horizon~$N = 30$ --~a property related to the degree of suboptimality~$\alpha_N$. For both the values of the prediction horizon, the asymptotic behavior of the PI-kEDMD-based MPC is the same that is obtained in a nominal MPC using the real system equations in the optimal control problem.
In Fig.~\ref{fig:vdp_comparison} (right), we report the trajectories of the same simulations in the phase space.
We can see that, while the asymptotic behavior is strongly related to the accuracy of the model in the origin, the transient behavior is mainly influenced by the prediction horizon and the number of cluster points.
Finally, Table~\ref{tab:times} reports the average computation time for the solution of the MPC on a laptop with an Intel i7-1065G7 CPU. The computation times increase when considering larger numbers of cluster point and longer prediction horizons, and are comparable for the PI-kEDMD and the kEDMD model.
\begin{table}[htb]
    \centering
    \bgroup
    \def\arraystretch{1.3}%
    \begin{tabular}{cc|c c | c c}
         & & \multicolumn{2}{c|}{PI-kEDMD} & \multicolumn{2}{c}{kEDMD} \\
         & & $N=10$ & $N=30$ & $N=10$ & $N=30$ \\
        \hline
        \multirow{2}{*}{\rotatebox{90}{V.d.P.}} &
        $d = 352$ & 0.039 & 0.114 & 0.049 & 0.118 \\
        & $d = 1327$ & 0.146 & 0.397 & 0.141 & 0.403 \\
        \hline
        \hline
        & & \multicolumn{2}{c|}{$N = 10$} & \multicolumn{2}{c}{$N = 10$} \\
        \hline
        \multirow{2}{*}{\rotatebox{90}{four-tanks\ }}
        & $d = 626$ & \multicolumn{2}{c|}{0.370} & \multicolumn{2}{c}{0.383}  \\
        & $d = 1297$ & \multicolumn{2}{c|}{0.875} & \multicolumn{2}{c}{0.864} \\
        & $d = 2402$ & \multicolumn{2}{c|}{1.876} & \multicolumn{2}{c}{1.744} 
    \end{tabular}
    \egroup
    \caption{Average computation time in seconds of the kEDMD-based MPC for the Van der Pol oscillator (V.d.P.) and for the four-tanks process, for different prediction horizons and numbers of data points.}
    \label{tab:times}
\end{table}

\subsection*{Four-tanks process}
Our second example is the four-tank system \cite{alvarado2011comparative} described by the differential equation
\begin{equation}\label{eq:fourTanks}
    \begin{pmatrix}
        \dot{h}_1 \\
        \dot{h}_2 \\
        \dot{h}_3 \\
        \dot{h}_4
    \end{pmatrix}
    = - \frac{\sqrt{2g}}{S} \begin{pmatrix}
        a_1 \sqrt{h_1} + a_3 \sqrt{h_3} \\
        a_2 \sqrt{h_2} + a_4 \sqrt{h_4} \\
        a_3 \sqrt{h_3} \\
        a_4 \sqrt{h_4}
    \end{pmatrix}
    + \begin{bmatrix}
        \frac{\gamma_a}{S} & 0 \\
        0 & \frac{\gamma_b}{S} \\       
        0 & \frac{1 - \gamma_b}{S} \\
        \frac{1 - \gamma_a}{S} & 0
    \end{bmatrix} \begin{pmatrix}
        q_a \\
        q_b
    \end{pmatrix}.
\end{equation}
The state of the system is given by the levels in the four tanks, i.e. $x = (h_1, \; h_2, \; h_3, \; h_4)^\top \in \R^4$, while the control variables are the flows through the two valves, i.e. $u = (q_a, \; q_b)^\top \in \R^2$. The numerical values of the system parameters can be found in \cite{alvarado2011comparative} and as sampling time we choose $\Delta t = 10s$. 
The discrete-time version of \eqref{eq:fourTanks} obtained with the forward Euler method is considered as ground truth. 
The controlled state for this system is the equilibrium given by $\bar{x} = (0.65, 0.66, 0.6417, 0.6882)^\top m$ and $\bar{u} = (1.666, 1.974)^\top m^3/h$. To derive the kEDMD surrogate model, the state and inputs are shifted around the origin. 
The state domain, where we sample, is $\Omega = [0.2\ m, 1.36\ m]^2 \times [0.2\ m, 1.30\ m]^2$, while the input is constrained by the set $\mathbb{U} = [0\ m^3/h, 3.26\ m^3/h]\times [0\ m^3/h, 4\ m^3/h]$. 
The cluster centers are chosen in a uniform grid in~$\Omega$. 
Moreover, an additional cluster center is added at the equilibrium. 
We consider datasets with $d \in \{626, 1297, 2402\}$ cluster centers, that correspond to $\sqrt[4]{d-1} \in \{5,6,7\}$. 
We decided to stick to a uniform grid, because the explicit construction of a grid minimizing the Lebesgue constant in dimensions higher than two is more involved and we did not run into numerical issues. 
Possible approaches to reduce the condition number of the model matrices are the use of the cartesian product of Padua grids or of more sophisticated approaches as outlined, e.g., in~\cite{jimenez2025approximating}.
The cluster radius is $\rx = 2/d$ for all the models with 25 samples $(x_{ij}, u_{ij}, x_{ij}^+)$ in each cluster.
For the closed-loop simulations, we consider an initial state $x_0 = (1.0, 1.0, 1.0, 1.0)^\top m$, MPC cost matrices $Q = I_4$ and $R = 10^{-4}I_2$, and a prediction horizon $N = 10$.
\begin{figure}[htb]
    \centering
    \includegraphics[width=0.5\linewidth]{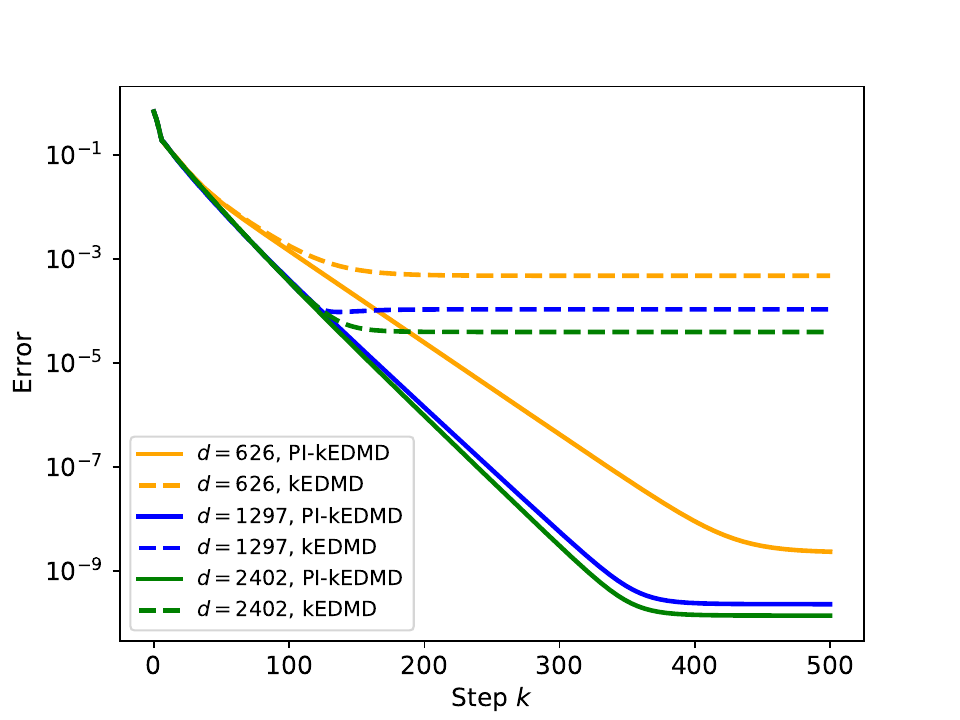}
    \caption{Four-tank: Error $\|x(k) - \bar{x}\|$ for the kEDMD-based MPC closed loop with horizon $N=10$.}
    \label{fig:four_tanks_N10}
\end{figure}

The errors of the closed loop simulations are reported in Fig. \ref{fig:four_tanks_N10}. 
As in the previous example, we can see that the error reaches a higher accuracy with the PI-kEDMD model, and the accuracy increases with a larger number of cluster points. The number of cluster points also influences how quickly the error decreases, which occurs more slowly for $d = 626$.
In this example, the closed loop behavior remains the same if the prediction horizon is increased to $N=20$. Finally, the average computation times for the solution of the MPC optimization problem are reported in Table~\ref{tab:times}.

\section{Conclusions and outlook}
\label{sec:conclusions}

\noindent We have proposed a framework to rigorously show asymptotic stability of an equilibrium w.r.t.\ the MPC closed, where data-driven surrogate models are used in the optimization step assuming cost controllability of the original system. 
The key ingredients are three properties of the data-driven surrogate: 
First, Lipschitz continuity uniform in the approximation accuracy. 
Second and third, uniform and foremost proportional error bounds, see Assumption~\ref{ass:error_bound}. 
Then, we rigorously verified all assumptions for a data-driven kEDMD surrogate models using Koopman operator theory resulting in finite-data error bounds. 
In particular, we have proposed a learning algorithm based on flexible data sampling. 
Finally, the theoretical results were verified by numerical simulations for two representative examples.

Future work could consider an extension towards Koopman-based control of partial differential equations~\cite{Deut24} as well as an extension towards novel (E)DMD variants such as eigenfunction-based approximation via (control) Liouville operators~\cite{RoseKama24}, residual DMD~\cite{ColbAyto23} or via stochastic Koopman approximants~\cite{WannMezi22} using recently proposed finite-data error bounds for kEDMD derived~\cite{HertPhil25}.
Moreover, the stability analysis can be extended to consider the presence of state constraints, disturbances and time-varying reference signals, which are not taken into account in this paper and require additional treatment as outlined in Remark~\ref{rem:state-constraints}. In addition, also an extension towards input-output data based on, e.g., \cite{ManzLimo20}, would be of interest.

\bibliographystyle{plain}        
\bibliography{Data-driven_MPC}

\end{document}